# Fourier series (based) multiscale method for computational analysis in science and engineering:

# VI. Fourier series multiscale solution for wave propagation in a beam with rectangular cross section

Weiming Sun[a*+] and Zimao Zhang[b]

**Abstract:** Fourier series multiscale method, a concise and efficient analytical approach for multiscale computation, will be developed out of this series of papers. In the sixth paper, exact analysis of the wave propagation in a beam with rectangular cross section is extended to a thorough multiscale analysis for a system of completely coupled second order linear differential equations for modal functions, where general boundary conditions are prescribed. For this purpose, the modal function each is expressed as a linear combination of the corner function, the two boundary functions and the internal function, to ensure the series expressions obtained uniformly convergent and termwise differentiable up to second order. Meanwhile, the sum of the corner function and the internal function corresponds to the particular solution, and the two boundary functions correspond to the general solutions which satisfy the homogeneous form of the equations. Since the general solutions have appropriately interpreted the meaning of the differential equations, the spatial characteristics of the solution of the equations are expected to be better captured in separate directions. With the corner function, the two boundary functions and the internal function selected specifically as polynomials, one-dimensional full-range Fourier series along the $x_2$ (or $x_1$)-direction, and two-dimensional full-range Fourier series, the Fourier series multiscale solution of the wave propagation in a beam with rectangular cross section is derived. For the beam with various boundary conditions, computation and analysis are performed, and the propagation characteristics of elastic waves in the beam are presented. The newly proposed accurate wave model has laid a solid foundation for simultaneous control of coupled waves in the beam and establishment of guided wave NDE techniques.

# 1. Introduction

The propagation-based active control concept has been put forward for almost 40 years,

[a]Department of Mathematics and Big Data, School of Artificial Intelligence, Jianghan University, Wuhan, 430056, China
[b]Department of Mechanics, School of Civil Engineering, Beijing Jiaotong University, Beijing, 100044, China
*Correspondence to: Weiming Sun, Department of Mathematics and Big Data, School of Artificial Intelligence, Jianghan University, Wuhan, 430056, China
[+] E-mail: xuxinenglish@hust.edu.cn

during which a number of works come out in this area. However, special stresses are all laid on such theoretical problems as design of control system, analysis of stability and choice of measuring methods of the fundamental continua [1-17]. By contrast, the analysis of wave propagation along the waveguides, despite its fundamental position in the traveling wave control approach, seems a little bit casual and oversimplified. For example, simultaneous control for different types of structural waves propagating in a beam has been studied by various researchers, where different types of waves were considered in isolation and controlled independently [18-20] at the expense of neglecting coupling effects and thus inevitably resulted in control spillover. Therefore, it can be seen that accurate determination of permissible waves within the waveguides is of great importance to the effective control of wave motion in structures. Meanwhile, ultrasonic nondestructive evaluation (NDE) has attracted considerable attention recently. In the application of guided waves in nondestructive evaluation, fundamental characteristics of guided waves, especially dispersion relations are often useful [21, 22]. Dispersion curves present fundamental information on guided waves such as wavelength and dispersivity as well as phase and group velocities at a certain frequency. And such fundamental information plays an important role in detection of defects in waveguides.

Problems associated with the guided wave propagation of harmonic waves in elastic wavegudes of various configurations have been the object of theoretical, experimental and numerical studies for more than 160 years [23-29]. This process might be described as a traditional dilemma "approximate solution of exact equations or exact solution of approximate equations". However, due to difficulties in mechanics and mathematics, up to now "exact solution of exact equations" can be carried out theoretically only in several cases, which include the infinite plate investigated by Rayleigh and Lamb in 1888, and the infinite beam of circular cross section dealt with, respectively, by Pochhammer in 1876 and Chree in 1889. As to the waveguides of other kinds, for instance, the infinite rectangular beam with stress-free boundaries, the similar exact treatment is considered impossible within the framework of linear elasticity theory [30, 31].

As a general analytical method for differential equations, the Fourier series method is of both great theoretical significance, and important practical value for the exact analysis of wave propagation in an infinite rectangular beam. For example, by solving the three-dimensional elastodynamic equations, Sun et al. derived the general analytical Fourier series solution for wave motion within the beam [32]. According to Lame's displacement potential [33], two kinds of waves can exist in an elastic medium: equivoluminal and rotationless waves, the linear combinations of which, as shown in the obtained series solution, also describe the wave motion in the beam. The series solution has rapid convergence speed and good accuracy not only for the first branches but also for the higher order branches of the dispersion curves. And then, this problem was investigated by Krushynska and Meleshko [34] in 2011. Small motions of the beam were expressed by Lame equations in terms of the components of displacement vector. Solution to these equations can be further obtained in terms of vector potential and scalar potential. By satisfying the decoupled Helmholtz equations for potentials, four independent families of normal waves (or homogenous solutions) were yielded. And the normal waves for each mode family were combined to obtain the solution for the corresponding wave mode in the beam. Numerical investigations were performed for various geometric and physical parameters of the beam, and a good agreement was found between their calculations and available numerical and experimental data.

Recently, on the basis of sufficient conditions of $2r$ ($r$ is a positive integer) times term-by-term differentiation of Fourier series [35], the methodology of simultaneous approximation with composite Fourier series has been proposed for (solution) functions and their (partial) derivatives [36]. In this methodology, a one-dimensional or two-dimensional (solution) function with general boundary conditions is decomposed into the linear combination of a corner function that describes the discontinuities at corners of the domain (only for the two-dimensional function), boundary functions that describe the discontinuities on boundaries of the domain and an internal function that describes the smoothness within the domain. Meanwhile, the obtained composite Fourier series possesses the flexibility of the formulation to a certain degree and allows of further adjustments in expression. Usually, the corner function and the internal function can be selected respectively as the algebraical polynomials with specific orders and Fourier series. However, the selection of interpolation basis functions for the boundary functions might differ in different situations. Specifically, in the process of solving the linear differential equation with constant coefficients, the selected interpolation basis functions satisfy the homogeneous form of the differential equation to be solved, and then the boundary functions constitute the general solutions of the original linear differential equation with constant coefficients. It can be expected that such general solutions shall be able to appropriately interpret the meaning of the differential equation, and hence better capture the spatial characteristics of the solution in separate directions. In view of the multiscale capability of this solution method, we rename it the Fourier series (based) multiscale method for linear differential equation with constant coefficients [37].

And now, a comparision of the structures of the Fourier series multiscale solution and the series solutions obtained in [32, 34] can be made. It is easy to observe that the series solutions obtained in [32, 34] virtually correspond to the general solutions (boundary functions) within the framework of the Fourier series multiscale method. Since discontinuities potentially related to the solutions and their partial derivatives on the boundaries (when they are periodically extended onto the entire $x_1 - x_2$ plane as implied by a two-dimensional Fourier series expansion) have not been explicited absorbed, there is the probability that the solutions for some specified boundary conditions cannot be differentiated term-by-term to obtain all the series expansions for up to the second order partial derivatives.

Therefore, in the sixth paper of the series of researches on Fourier series multiscale method, we will reinvestigate the wave propagation in a beam with rectangular cross section from a more systematic and general view of point. Firstly, we generalize this issue to the multiscale analysis of a system of completely coupled second order linear differential equations (for modal functions), where general boundary conditions are prescribed. Secondly, the Fourier series multiscale method is extended for the exact analysis of the wave propagation in a beam with rectangular cross section. As a routine task, the modal functions are respectively decomposed into several constituents, such as the corner function, the two boundary functions and the internal function. The sum of the corner function and the internal function corresponds to the particular solutions. And the two boundary functions correspond to the general solutions which satisfy the homogeneous form of the governing differential equation. And then, with all the constituents selected respectively as polynomials, one-dimensional full-range Fourier series along the $x_2$ (or $x_1$)-direction, and two-dimensional full-range Fourier series, the Fourier series multiscale solution of wave propagation in a beam with rectangular cross section is derived. Accordingly, this paper begins with the three-dimensional elastodynamic equations. Detailed formulations related to the Fourier series multiscale solution are then presented. Finally, computation and analysis are performed for the beam under various boundary conditions, where the propagation characteristics of the frequency spectrum, cut-off frequencies and wave modes are given. The newly proposed accurate wave model has laid a solid foundation for simultaneous control of coupled waves in the beam and establishment of guided wave NDE techniques.

## 2. Description of the problem

As shown in Figure 1, an infinite beam with Lame's constants of $\lambda$, $\mu$ and rectangular cross section of $2a$ in length, $2b$ in width is considered.

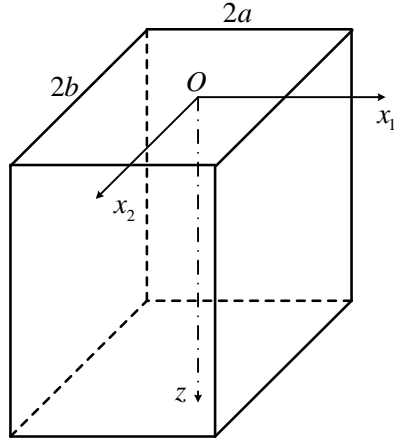

Figure 1: An infinite beam of rectangular cross section.

Accordingly, the three-dimensional elastodynamic equation, in the form of displacements, is given below [30]

$$\left.\begin{aligned}\rho\frac{\partial^2 u}{\partial t^2} &= (\lambda+\mu)\frac{\partial e}{\partial x_1}+\mu\nabla^2 u \\ \rho\frac{\partial^2 v}{\partial t^2} &= (\lambda+\mu)\frac{\partial e}{\partial x_2}+\mu\nabla^2 v \\ \rho\frac{\partial^2 w}{\partial t^2} &= (\lambda+\mu)\frac{\partial e}{\partial z}+\mu\nabla^2 w\end{aligned}\right\}, \qquad (1)$$

where $u$, $v$ and $w$ are the displacement components at any point of the elastic body, $\rho$ is the material density, $\omega$ is the harmonic frequency, and the volumetric strain

$$e = \frac{\partial u}{\partial x_1}+\frac{\partial v}{\partial x_2}+\frac{\partial w}{\partial z}, \qquad (2)$$

and the Laplace differential operator

$$\nabla^2 = \frac{\partial^2}{\partial x_1^2}+\frac{\partial^2}{\partial x_2^2}+\frac{\partial^2}{\partial z^2}. \qquad (3)$$

The corresponding constitutive equation is

$$\left.\begin{aligned}\sigma_{x_1} &= \lambda e+2\mu\frac{\partial u}{\partial x_1}, \quad \sigma_{x_2} = \lambda e+2\mu\frac{\partial v}{\partial x_2}, \quad \sigma_z = \lambda e+2\mu\frac{\partial w}{\partial z} \\ \tau_{x_1 x_2} &= \mu(\frac{\partial u}{\partial x_2}+\frac{\partial v}{\partial x_1}), \quad \tau_{x_2 z} = \mu(\frac{\partial v}{\partial z}+\frac{\partial w}{\partial x_2}), \quad \tau_{x_1 z} = \mu(\frac{\partial w}{\partial x_1}+\frac{\partial u}{\partial z})\end{aligned}\right\}. \qquad (4)$$

For example, there are two kinds of typical boundary conditions at the edge $x_1 = a$ of the beam:

1. Clamped edge (C for short)

$$u(a, x_2) = 0, \quad v(a, x_2) = 0, \quad w(a, x_2) = 0, \qquad (5)$$

2. Free edge (F for short)

$$\sigma_{x_1}(a, x_2) = 0, \quad \tau_{x_1 x_2}(a, x_2) = 0, \quad \tau_{x_1 z}(a, x_2) = 0. \tag{6}$$

## 3. The Fourier series multiscale solution

The wave propagation in a rectangular beam can be formulated as a system of differential equations consisting of second order ($2r_u = 2r_v = 2r_w = 2$) linear differential equations with constant coefficients of the spatial displacement components $u$, $v$ and $w$. In addition, the system of differential equations contains the solution functions and their mixed second order partial derivatives. Therefore, we can obtain the corresponding Fourier series multiscale solution in the form of the two-dimensional full-range expansion [37] of the spatial displacement components.

*3.1. System of differential equations of modal functions*

We now confine ourselves to the harmonic wave, which propagates along the beam without changing its form. Therefore, suppose that Eq. (1) has such displacement solution as

$$\left.\begin{aligned} u &= \varphi_u(x_1, x_2) \cos(kz - \omega t) \\ v &= \varphi_v(x_1, x_2) \cos(kz - \omega t) \\ w &= \varphi_w(x_1, x_2) \sin(kz - \omega t) \end{aligned}\right\}, \tag{7}$$

where $\varphi_u$, $\varphi_v$ and $\varphi_w$ are usually termed the modal functions and the wave number $k$ might be real, imaginary and complex, and these three cases correspond respectively to the traveling, near field and attenuate traveling wave in the beam.

We then substitute the displacement solution into Eq. (1), and obtain the system of differential equations of the modal functions

$$\mathbf{L}_\varphi \begin{bmatrix} \varphi_u \\ \varphi_v \\ \varphi_w \end{bmatrix} = 0, \tag{8}$$

where the matrix of differential operators

$$\mathbf{L}_\varphi = \begin{bmatrix} \mathbf{L}_{\varphi u} \\ \mathbf{L}_{\varphi v} \\ \mathbf{L}_{\varphi w} \end{bmatrix}, \tag{9}$$

and the vectors of differential operators

$$\mathbf{L}_{\varphi u} = \left[ (\lambda + 2\mu)\frac{\partial^2}{\partial x_1^2} + \mu\frac{\partial^2}{\partial x_2^2} + (\rho\omega^2 - \mu k^2) \quad (\lambda + \mu)\frac{\partial^2}{\partial x_1 \partial x_2} \quad (\lambda + \mu)k\frac{\partial}{\partial x_1} \right], \tag{10}$$

$$\mathbf{L}_{\varphi v} = \left[ (\lambda + \mu)\frac{\partial^2}{\partial x_1 \partial x_2} \quad \mu\frac{\partial^2}{\partial x_1^2} + (\lambda + 2\mu)\frac{\partial^2}{\partial x_2^2} + (\rho\omega^2 - \mu k^2) \quad (\lambda + \mu)k\frac{\partial}{\partial x_2} \right], \tag{11}$$

$$\mathbf{L}_{\varphi w} = \left[ -(\lambda + \mu)k\frac{\partial}{\partial x_1} \quad -(\lambda + \mu)k\frac{\partial}{\partial x_2} \quad \mu\frac{\partial^2}{\partial x_1^2} + \mu\frac{\partial^2}{\partial x_2^2} + \rho\omega^2 - (\lambda + 2\mu)k^2 \right]. \tag{12}$$

## 3.2. Structural decomposition of modal functions

According to the Fourier series multiscale method for the second order linear differential equation with constant coefficients [37], we expand the modal functions $\varphi_u$, $\varphi_v$ and $\varphi_w$ in composite full-range Fourier series over the domain $[-a,a] \times [-b,b]$

$$\varphi_u(x_1, x_2) = \varphi_{0u}(x_1, x_2) + \varphi_{1u}(x_1, x_2) + \varphi_{2u}(x_1, x_2) + \varphi_{3u}(x_1, x_2)$$
$$= \mathbf{\Phi}_0^T(x_1, x_2) \cdot \mathbf{q}_{0,u} + \mathbf{\Phi}_{1,u}^T(x_1, x_2) \cdot \mathbf{q}_{1,u} + \mathbf{\Phi}_{2,u}^T(x_1, x_2) \cdot \mathbf{q}_{2,u} + \mathbf{\Phi}_3^T(x_1, x_2) \cdot \mathbf{q}_{3,u}, \quad (13)$$

$$\varphi_v(x_1, x_2) = \varphi_{0v}(x_1, x_2) + \varphi_{1v}(x_1, x_2) + \varphi_{2v}(x_1, x_2) + \varphi_{3v}(x_1, x_2)$$
$$= \mathbf{\Phi}_0^T(x_1, x_2) \cdot \mathbf{q}_{0,v} + \mathbf{\Phi}_{1,v}^T(x_1, x_2) \cdot \mathbf{q}_{1,v} + \mathbf{\Phi}_{2,v}^T(x_1, x_2) \cdot \mathbf{q}_{2,v} + \mathbf{\Phi}_3^T(x_1, x_2) \cdot \mathbf{q}_{3,v}, \quad (14)$$

$$\varphi_w(x_1, x_2) = \varphi_{0w}(x_1, x_2) + \varphi_{1w}(x_1, x_2) + \varphi_{2w}(x_1, x_2) + \varphi_{3w}(x_1, x_2)$$
$$= \mathbf{\Phi}_0^T(x_1, x_2) \cdot \mathbf{q}_{0,w} + \mathbf{\Phi}_{1,w}^T(x_1, x_2) \cdot \mathbf{q}_{1,w} + \mathbf{\Phi}_{2,w}^T(x_1, x_2) \cdot \mathbf{q}_{2,w} + \mathbf{\Phi}_3^T(x_1, x_2) \cdot \mathbf{q}_{3,w}, \quad (15)$$

where the vectors of basis functions, i.e., $\mathbf{\Phi}_0^T$, $\mathbf{\Phi}_3^T$, and the undetermined constant vectors, i.e., $\mathbf{q}_{i,u}$, $\mathbf{q}_{i,v}$, $\mathbf{q}_{i,w}$, $i=0,1,2,3$, are as defined in section 3.2, [36]. Moreover, the specific expressions of the vectors of basis functions $\mathbf{\Phi}_{1,u}^T$, $\mathbf{\Phi}_{2,u}^T$, $\mathbf{\Phi}_{1,v}^T$, $\mathbf{\Phi}_{2,v}^T$, $\mathbf{\Phi}_{1,w}^T$ and $\mathbf{\Phi}_{2,w}^T$ are to be determined with the general solution of Eq. (8).

## 3.3. Expressions of boundary functions expanded along the $x_2$-direction

Suppose that Eq. (8) has the following homogeneous solution

$$\left.\begin{array}{l}\varphi_{1n,u}(x_1, x_2) = \xi_{1n,u}(x_1)\cos(\beta_n x_2) + \xi_{2n,u}(x_1)\sin(\beta_n x_2) \\ \varphi_{1n,v}(x_1, x_2) = \xi_{1n,v}(x_1)\cos(\beta_n x_2) + \xi_{2n,v}(x_1)\sin(\beta_n x_2) \\ \varphi_{1n,w}(x_1, x_2) = \xi_{1n,w}(x_1)\cos(\beta_n x_2) + \xi_{2n,w}(x_1)\sin(\beta_n x_2)\end{array}\right\}, \quad (16)$$

where $n$ is a nonnegative integer, the functions $\xi_{1n,u}(x_1)$, $\xi_{2n,u}(x_1)$, $\xi_{1n,v}(x_1)$, $\xi_{2n,v}(x_1)$, $\xi_{1n,w}(x_1)$ and $\xi_{2n,w}(x_1)$ are undetermined one-dimensional functions, and $\beta_n = n\pi/b$.

Substituting Eq. (16) into Eq. (8), we then obtain that

1. for $n > 0$, the system of equations for the undetermined functions $\xi_{1n,u}(x_1)$, $\xi_{2n,v}(x_1)$ and $\xi_{1n,w}(x_1)$

$$\mathbf{L}_{1n}\begin{bmatrix}\xi_{1n,u} \\ \xi_{2n,v} \\ \xi_{1n,w}\end{bmatrix} = 0, \quad (17)$$

where the matrix of differential operators

$$\mathbf{L}_{1n} = \begin{bmatrix}\mathbf{L}_{1n,u} \\ \mathbf{L}_{1n,v} \\ \mathbf{L}_{1n,w}\end{bmatrix}, \quad (18)$$

and the vectors of differential operators

$$\mathbf{L}_{1n,u} = \left[(\lambda + 2\mu)\frac{d^2}{dx_1^2} + \rho\omega^2 - \mu\beta_n^2 - \mu k^2 \quad (\lambda+\mu)\beta_n\frac{d}{dx_1} \quad (\lambda+\mu)k\frac{d}{dx_1}\right], \quad (19)$$

$$\mathbf{L}_{1n,v} = \left[ -(\lambda+\mu)\beta_n \frac{d}{dx_1} \quad \mu\frac{d^2}{dx_1^2} + \rho\omega^2 - (\lambda+2\mu)\beta_n^2 - \mu k^2 \quad -(\lambda+\mu)\beta_n k \right], \quad (20)$$

$$\mathbf{L}_{1n,w} = \left[ -(\lambda+\mu)k \frac{d}{dx_1} \quad -(\lambda+\mu)\beta_n k \quad \mu\frac{d^2}{dx_1^2} + \rho\omega^2 - \mu\beta_n^2 - (\lambda+2\mu)k^2 \right]. \quad (21)$$

2. for $n > 0$, the system of equations for the undetermined functions $\xi_{2n,u}(x_1)$, $\xi_{1n,v}(x_1)$ and $\xi_{2n,w}(x_1)$

$$\mathbf{L}_{2n} \begin{bmatrix} \xi_{2n,u} \\ \xi_{1n,v} \\ \xi_{2n,w} \end{bmatrix} = 0, \quad (22)$$

where the matrix of differential operators $\mathbf{L}_{2n}$ has similar form as that of $\mathbf{L}_{1n}$, with the parameter $\beta_n$ substituted by $-\beta_n$.

3. for $n = 0$, the system of equations for the undetermined functions $\xi_{10,u}(x_1)$, $\xi_{10,v}(x_1)$ and $\xi_{10,w}(x_1)$

$$\mathbf{L}_{10} \begin{bmatrix} \xi_{10,u} \\ \xi_{10,v} \\ \xi_{10,w} \end{bmatrix} = 0, \quad (23)$$

where the matrix of differential operators $\mathbf{L}_{10}$ has similar form as that of $\mathbf{L}_{1n}$, but to substitute the parameter $\beta_n$ in $\mathbf{L}_{1n}$ with 0.

1. Solving Eq. (17)

Further, let the exponential solution of Eq. (17) be

$$\left. \begin{array}{l} \xi_{1n,u}(x_1) = G_{1n1,u} p_{1n}(x_1) \\ \xi_{2n,u}(x_1) = G_{1n2,v} p_{1n}(x_1) \\ \xi_{1n,w}(x_1) = G_{1n1,w} p_{1n}(x_1) \end{array} \right\}, \quad (24)$$

where the basis function

$$p_{1n}(x_1) = \exp(\eta_n x_1), \quad (25)$$

$G_{1n1,u}$, $G_{1n2,v}$, $G_{1n1,w}$ and $\eta_n$ are undetermined constants.

Substituting Eq. (25) in Eq. (17), we obtain the characteristics equation

$$[\mu(-\eta_n^2 + \beta_n^2 + k^2) - \rho w^2]^2 [(\lambda+2\mu)(-\eta_n^2 + \beta_n^2 + k^2) - \rho w^2] = 0. \quad (26)$$

Therefore, it is easy to verify that Eq. (26) has two double roots and two single roots. Particularly, if $k$ is a real number, we denote the discriminant as

$$\Delta_1 = \beta_n^2 + k^2 - \frac{\rho w^2}{\mu}. \quad (27)$$

Then if $\Delta_1 > 0$, the two double roots of Eq. (26) can be expressed as

$$\eta_{n,1} = \eta_{n,2} = \alpha_{1n}, \quad \eta_{n,3} = \eta_{n,4} = -\alpha_{1n}, \quad (28)$$

where $\alpha_{1n} = \sqrt{\beta_n^2 + k^2 - \frac{\rho w^2}{\mu}}$.

If $\Delta_1 < 0$, the two double roots of Eq. (26) can be expressed as

$$\eta_{n,1} = \eta_{n,2} = i\alpha_{2n}, \quad \eta_{n,3} = \eta_{n,4} = -i\alpha_{2n}, \tag{29}$$

where $\alpha_{2n} = \sqrt{-(\beta_n^2 + k^2 - \dfrac{\rho w^2}{\mu})}$.

In addition, we denote the discriminant as

$$\Delta_2 = \beta_n^2 + k^2 - \frac{\rho w^2}{\lambda + 2\mu}. \tag{30}$$

Then, if $\Delta_2 > 0$, the two single roots of Eq. (26) can be expressed as

$$\eta_{n,5} = \alpha_{3n}, \quad \eta_{n,6} = -\alpha_{3n}, \tag{31}$$

where $\alpha_{3n} = \sqrt{\beta_n^2 + k^2 - \dfrac{\rho w^2}{\lambda + 2\mu}}$.

If $\Delta_2 < 0$, the two single roots of Eq. (26) can be expressed as

$$\eta_{n,5} = i\alpha_{4n}, \quad \eta_{n,6} = -i\alpha_{4n}, \tag{32}$$

where $\alpha_{4n} = \sqrt{-(\beta_n^2 + k^2 - \dfrac{\rho w^2}{\lambda + 2\mu})}$.

And accordingly, the basis functions $p_{1nl}(x_1)$, $l = 1, 2, 3, 4$, are presented in Table 1.

Table 1: Expressions for $p_{1nl}(x_1)$, $l = 1, 2, 3, 4$.

|  | $\Delta_1 > 0$ | $\Delta_1 < 0$ | $\Delta_2 > 0$ | $\Delta_2 < 0$ |
|---|---|---|---|---|
| $p_{1n1}(x_1)$ | $\cosh(\alpha_{1n}x_1)$ | $\cos(\alpha_{2n}x_1)$ |  |  |
| $p_{1n2}(x_1)$ | $\sinh(\alpha_{1n}x_1)$ | $\sin(\alpha_{2n}x_1)$ |  |  |
| $p_{1n3}(x_1)$ |  |  | $\cosh(\alpha_{3n}x_1)$ | $\cos(\alpha_{4n}x_1)$ |
| $p_{1n4}(x_1)$ |  |  | $\sinh(\alpha_{3n}x_1)$ | $\sin(\alpha_{4n}x_1)$ |

With the same procedure in Section 4.2, [37], we select the following vector of functions

$$\mathbf{p}_{1n}^{\mathrm{T}}(x_1) = [p_{1n1}(x_1) \quad p_{1n2}(x_1) \quad p_{1n3}(x_1) \quad p_{1n4}(x_1)], \tag{33}$$

and further, let the undetermined constant vectors be

$$\mathbf{a}_{1n1,u}^{\mathrm{T}} = [G_{1n1,u}^1 \quad G_{1n1,u}^2 \quad G_{1n1,u}^3 \quad G_{1n1,u}^4], \tag{34}$$

$$\mathbf{a}_{1n2,v}^{\mathrm{T}} = [G_{1n2,v}^1 \quad G_{1n2,v}^2 \quad G_{1n2,v}^3 \quad G_{1n2,v}^4], \tag{35}$$

$$\mathbf{a}_{1n1,w}^{\mathrm{T}} = [G_{1n1,w}^1 \quad G_{1n1,w}^2 \quad G_{1n1,w}^3 \quad G_{1n1,w}^4]. \tag{36}$$

Thus, we can construct functions $\xi_{1n,u}(x_1)$, $\xi_{2n,v}(x_1)$ and $\xi_{1n,w}(x_1)$ as follows

$$\xi_{1n,u}(x_1) = \mathbf{p}_{1n}^{\mathrm{T}}(x_1) \cdot \mathbf{a}_{1n1,u}, \tag{37}$$

$$\xi_{2n,v}(x_1) = \mathbf{p}_{1n}^{\mathrm{T}}(x_1) \cdot \mathbf{a}_{1n2,v}, \tag{38}$$

$$\xi_{1n,w}(x_1) = \mathbf{p}_{1n}^{\mathrm{T}}(x_1) \cdot \mathbf{a}_{1n1,w}. \tag{39}$$

Substituting the above expressions in Eq. (17), we obtain the following relations of the undetermined constant vectors $\mathbf{a}_{1n1,u}^{\mathrm{T}}$, $\mathbf{a}_{1n2,v}^{\mathrm{T}}$ and $\mathbf{a}_{1n1,w}^{\mathrm{T}}$:

a. if $\Delta_1 > 0$, then

$$\left.\begin{array}{l}\alpha_{1n}G^1_{1n1,u} + \beta_n G^2_{1n2,v} + kG^2_{1n1,w} = 0 \\ \alpha_{1n}G^2_{1n1,u} + \beta_n G^1_{1n2,v} + kG^1_{1n1,w} = 0\end{array}\right\}, \tag{40}$$

b. if $\Delta_1 < 0$, then

$$\left.\begin{array}{l}-\alpha_{2n}G^1_{1n1,u} + \beta_n G^2_{1n2,v} + kG^2_{1n1,w} = 0 \\ \alpha_{2n}G^2_{1n1,u} + \beta_n G^1_{1n2,v} + kG^1_{1n1,w} = 0\end{array}\right\}, \tag{41}$$

c. if $\Delta_2 > 0$, then

$$\left.\begin{array}{l}\dfrac{G^3_{1n1,u}}{-\alpha_{3n}} = \dfrac{G^4_{1n2,v}}{\beta_n} = \dfrac{G^4_{1n1,w}}{k} \\ \dfrac{G^4_{1n1,u}}{-\alpha_{3n}} = \dfrac{G^3_{1n2,v}}{\beta_n} = \dfrac{G^3_{1n1,w}}{k}\end{array}\right\}, \tag{42}$$

d. if $\Delta_2 < 0$, then

$$\left.\begin{array}{l}\dfrac{G^3_{1n1,u}}{-\alpha_{4n}} = \dfrac{G^4_{1n2,v}}{\beta_n} = \dfrac{G^4_{1n1,w}}{k} \\ \dfrac{G^4_{1n1,u}}{\alpha_{4n}} = \dfrac{G^3_{1n2,v}}{\beta_n} = \dfrac{G^3_{1n1,w}}{k}\end{array}\right\}. \tag{43}$$

If we denote the undetermined constant vectors as

$$\mathbf{a}^T_{1n1} = [\mathbf{a}^T_{1n1,u} \quad \mathbf{a}^T_{1n2,v} \quad \mathbf{a}^T_{1n1,w}], \tag{44}$$

$$\mathbf{a}^T_{1n1,R} = [G^1_{1n1,u} \quad G^2_{1n1,u} \quad G^3_{1n2,v} \quad G^4_{1n2,v} \quad G^1_{1n1,w} \quad G^2_{1n1,w}], \tag{45}$$

then by combining Eqs. (40)-(43), we have

$$\mathbf{a}_{1n1} = \mathbf{T}_{1n1}\mathbf{a}_{1n1,R}, \tag{46}$$

where $\mathbf{T}_{1n1}$ is a transformation matrix.

In addition, referring to Eqs. (62) and (63) in [36], we obtain

$$\mathbf{R}_{1n}\mathbf{a}_{1n1,u} = \mathbf{q}_{1n1,u}, \tag{47}$$

$$\mathbf{R}_{1n}\mathbf{a}_{1n2,v} = \mathbf{q}_{1n2,v}, \tag{48}$$

$$\mathbf{R}_{1n}\mathbf{a}_{1n1,u} = \mathbf{q}_{1n1,u}, \tag{49}$$

where the matrix

$$\mathbf{R}_{1n} = \begin{bmatrix}\mathbf{p}^{(0)T}_{1n}(a) - \mathbf{p}^{(0)T}_{1n}(-a) \\ \mathbf{p}^{(1)T}_{1n}(a) - \mathbf{p}^{(1)T}_{1n}(-a)\end{bmatrix}, \tag{50}$$

$\mathbf{q}_{1n1,u}$, $\mathbf{q}_{1n2,v}$ and $\mathbf{q}_{1n1,w}$ are the boundary Fourier coefficient subvectors, and are as defined in Eqs. (65) and (66), [36].

Further, Eqs. (47)-(49) can be rewritten as

$$\mathbf{S}_{1n}\mathbf{a}_{1n1} = \mathbf{q}_{1n1}, \tag{51}$$

where the matrix

$$\mathbf{S}_{1n} = \begin{bmatrix}\mathbf{R}_{1n} & \mathbf{0} & \mathbf{0} \\ \mathbf{0} & \mathbf{R}_{1n} & \mathbf{0} \\ \mathbf{0} & \mathbf{0} & \mathbf{R}_{1n}\end{bmatrix}, \tag{52}$$

and the boundary Fourier coefficient subvector
$$\mathbf{q}_{1n1}^{T} = [\mathbf{q}_{1n1,u}^{T} \quad \mathbf{q}_{1n2,v}^{T} \quad \mathbf{q}_{1n1,w}^{T}]. \tag{53}$$

Combining Eqs. (46) and (51), we have
$$\mathbf{S}_{1n}\mathbf{T}_{1n1}\mathbf{a}_{1n1,R} = \mathbf{q}_{1n1}, \tag{54}$$

then
$$\mathbf{a}_{1n1,R} = (\mathbf{S}_{1n}\mathbf{T}_{1n1})^{-1}\mathbf{q}_{1n1}. \tag{55}$$

Therefore, the undetermined one-dimensional functions $\xi_{1n,u}(x_1)$, $\xi_{2n,v}(x_1)$ and $\xi_{1n,w}(x_1)$ can eventually be expressed as

$$\begin{bmatrix} \xi_{1n,u} \\ \xi_{2n,v} \\ \xi_{1n,w} \end{bmatrix} = \mathbf{p}_{1n1,R}^{T}(x_1) \cdot \mathbf{q}_{1n1}, \tag{56}$$

where the matrix
$$\mathbf{p}_{1n1,R}^{T}(x_1) = \begin{bmatrix} \mathbf{p}_{1n}^{T}(x_1) & \mathbf{0} & \mathbf{0} \\ \mathbf{0} & \mathbf{p}_{1n}^{T}(x_1) & \mathbf{0} \\ \mathbf{0} & \mathbf{0} & \mathbf{p}_{1n}^{T}(x_1) \end{bmatrix} \cdot \mathbf{T}_{1n1} \cdot (\mathbf{S}_{1n}\mathbf{T}_{1n1})^{-1}. \tag{57}$$

Accordingly, we define the matrix of basis functions as
$$\mathbf{\Phi}_{1n,1}^{T}(x_1, x_2) = \mathbf{H}_{1n,1}(x_2) \cdot \mathbf{p}_{1n1,R}^{T}(x_1), \tag{58}$$

where the matrix
$$\mathbf{H}_{1n,1}(x_2) = \begin{bmatrix} \cos(\beta_n x_2) & 0 & 0 \\ 0 & \sin(\beta_n x_2) & 0 \\ 0 & 0 & \cos(\beta_n x_2) \end{bmatrix}, \tag{59}$$

and the corresponding matrix of higher order partial derivatives of the basis functions as
$$\mathbf{\Phi}_{1n,1}^{(k_1,k_2)T}(x_1, x_2) = \mathbf{H}_{1n,1}^{(k_2)}(x_2) \cdot \mathbf{p}_{1n1,R}^{(k_1)T}(x_1), \tag{60}$$

where $k_1$ and $k_2$ are nonnegative integers, and the matrix
$$\mathbf{p}_{1n1,R}^{(k_1)T}(x_1) = \begin{bmatrix} \mathbf{p}_{1n}^{(k_1)T}(x_1) & \mathbf{0} & \mathbf{0} \\ \mathbf{0} & \mathbf{p}_{1n}^{(k_1)T}(x_1) & \mathbf{0} \\ \mathbf{0} & \mathbf{0} & \mathbf{p}_{1n}^{(k_1)T}(x_1) \end{bmatrix} \cdot \mathbf{T}_{1n1} \cdot (\mathbf{S}_{1n}\mathbf{T}_{1n1})^{-1}, \tag{61}$$

and the matrix
$$\mathbf{H}_{1n,1}^{(k_2)}(x_2) = \begin{bmatrix} [\cos(\beta_n x_2)]^{(k_2)} & 0 & 0 \\ 0 & [\sin(\beta_n x_2)]^{(k_2)} & 0 \\ 0 & 0 & [\cos(\beta_n x_2)]^{(k_2)} \end{bmatrix}. \tag{62}$$

2. Solving Eq. (22)

Similarly, the undetermined one-dimensional functions $\xi_{2n,u}(x_1)$, $\xi_{1n,v}(x_1)$ and $\xi_{2n,w}(x_1)$ can eventually be expressed as

$$\begin{bmatrix} \xi_{2n,u} \\ \xi_{1n,v} \\ \xi_{2n,w} \end{bmatrix} = \mathbf{p}_{1n2,R}^{T}(x_1) \cdot \mathbf{q}_{1n2}, \tag{63}$$

where the matrix

$$\mathbf{p}_{1n2,R}^{\mathrm{T}}(x_1) = \begin{bmatrix} \mathbf{p}_{1n}^{\mathrm{T}}(x_1) & \mathbf{0} & \mathbf{0} \\ \mathbf{0} & \mathbf{p}_{1n}^{\mathrm{T}}(x_1) & \mathbf{0} \\ \mathbf{0} & \mathbf{0} & \mathbf{p}_{1n}^{\mathrm{T}}(x_1) \end{bmatrix} \cdot \mathbf{T}_{1n2} \cdot (\mathbf{S}_{1n}\mathbf{T}_{1n2})^{-1}, \tag{64}$$

the transformation matrix $\mathbf{T}_{1n2}$ has similar form as that of $\mathbf{T}_{1n1}$, with parameter $\beta_n$ substituted by $-\beta_n$ in $\mathbf{T}_{1n1}$, and the boundary Fourier coefficient subvector is

$$\mathbf{q}_{1n2}^{\mathrm{T}} = [\mathbf{q}_{1n2,u}^{\mathrm{T}} \quad \mathbf{q}_{1n1,v}^{\mathrm{T}} \quad \mathbf{q}_{1n2,w}^{\mathrm{T}}]. \tag{65}$$

Accordingly, we define the vector of basis functions as

$$\mathbf{\Phi}_{1n,2}^{\mathrm{T}}(x_1, x_2) = \mathbf{H}_{1n,2}(x_2) \cdot \mathbf{p}_{1n2,R}^{\mathrm{T}}(x_1), \tag{66}$$

where the function matrix

$$\mathbf{H}_{1n,2}(x_2) = \begin{bmatrix} \sin(\beta_n x_2) & 0 & 0 \\ 0 & \cos(\beta_n x_2) & 0 \\ 0 & 0 & \sin(\beta_n x_2) \end{bmatrix}, \tag{67}$$

and the matrix of higher order partial derivatives of the basis functions as

$$\mathbf{\Phi}_{1n,2}^{(k_1,k_2)\mathrm{T}}(x_1, x_2) = \mathbf{H}_{1n,2}^{(k_2)}(x_2) \cdot \mathbf{p}_{1n2,R}^{(k_1)\mathrm{T}}(x_1), \tag{68}$$

where $k_1$ and $k_2$ are nonnegative integers, and the matrixes $\mathbf{p}_{1n2,R}^{(k_1)\mathrm{T}}(x_1)$ and $\mathbf{H}_{1n,2}^{(k_2)}(x_2)$ are similar to those of Eqs. (61) and (62).

3. Solving Eq. (23)

Let the undetermined constant vectors be

$$\mathbf{a}_{101,u}^{\mathrm{T}} = [G_{101,u}^1 \quad G_{101,u}^2 \quad G_{101,u}^3 \quad G_{101,u}^4], \tag{69}$$

$$\mathbf{a}_{101,v}^{\mathrm{T}} = [G_{101,v}^1 \quad G_{101,v}^2 \quad G_{101,v}^3 \quad G_{101,v}^4], \tag{70}$$

$$\mathbf{a}_{101,w}^{\mathrm{T}} = [G_{101,w}^1 \quad G_{101,w}^2 \quad G_{101,w}^3 \quad G_{101,w}^4]. \tag{71}$$

We further establish the relations of the undetermined constant vectors $\mathbf{a}_{101,u}^{\mathrm{T}}$, $\mathbf{a}_{101,v}^{\mathrm{T}}$, $\mathbf{a}_{101,w}^{\mathrm{T}}$ as follows:

a. if $\Delta_1 > 0$, then

$$\left. \begin{array}{l} \alpha_{10} G_{101,u}^1 + k G_{101,w}^2 = 0 \\ \alpha_{10} G_{101,u}^2 + k G_{101,w}^1 = 0 \end{array} \right\}, \tag{72}$$

b. if $\Delta_1 < 0$, then

$$\left. \begin{array}{l} -\alpha_{20} G_{101,u}^1 + k G_{101,w}^2 = 0 \\ \alpha_{20} G_{101,u}^2 + k G_{101,w}^1 = 0 \end{array} \right\}, \tag{73}$$

c. if $\Delta_2 > 0$, then

$$\left. \begin{array}{l} \dfrac{G_{101,u}^3}{-\alpha_{30}} = \dfrac{G_{101,w}^4}{k} \\ \dfrac{G_{101,u}^4}{-\alpha_{30}} = \dfrac{G_{101,w}^3}{k} \end{array} \right\}, \tag{74}$$

d. if $\Delta_2 < 0$, then

$$\left.\begin{array}{l}\dfrac{G_{101,u}^{3}}{-\alpha_{40}}=\dfrac{G_{101,w}^{4}}{k}\\[6pt]\dfrac{G_{101,u}^{4}}{\alpha_{40}}=\dfrac{G_{101,w}^{3}}{k}\end{array}\right\}, \tag{75}$$

e. $G_{101,v}^{1}$ and $G_{101,v}^{2}$ are arbitrary constants, (76)

f. $G_{101,v}^{3}=G_{101,v}^{4}=0$. (77)

If we define the undetermined constant vectors as

$$\mathbf{a}_{101}^{T}=[\mathbf{a}_{101,u}^{T}\quad \mathbf{a}_{101,v}^{T}\quad \mathbf{a}_{101,w}^{T}], \tag{78}$$

$$\mathbf{a}_{101,R}^{T}=[G_{101,u}^{1}\quad G_{101,u}^{2}\quad G_{101,v}^{1}\quad G_{101,v}^{2}\quad G_{101,w}^{3}\quad G_{101,w}^{4}], \tag{79}$$

then Eqs. (72)-(77) are combined to be

$$\mathbf{a}_{101}=\mathbf{T}_{101}\mathbf{a}_{101,R}, \tag{80}$$

where $\mathbf{T}_{101}$ is a transformation matrix.

In addition, the undetermined one-dimensional functions $\xi_{10,u}(x_1)$, $\xi_{10,v}(x_1)$ and $\xi_{10,w}(x_1)$ are expressed as

$$\begin{bmatrix}\xi_{10,u}\\ \xi_{10,v}\\ \xi_{10,w}\end{bmatrix}=\mathbf{p}_{101,R}^{T}(x_1)\cdot\mathbf{q}_{101}, \tag{81}$$

where the matrix

$$\mathbf{p}_{101,R}^{T}(x_1)=\begin{bmatrix}\mathbf{p}_{10}^{T}(x_1) & \mathbf{0} & \mathbf{0}\\ \mathbf{0} & \mathbf{p}_{10}^{T}(x_1) & \mathbf{0}\\ \mathbf{0} & \mathbf{0} & \mathbf{p}_{10}^{T}(x_1)\end{bmatrix}\cdot\mathbf{T}_{101}\cdot(\mathbf{S}_{10}\mathbf{T}_{101})^{-1}, \tag{82}$$

and the boundary Fourier coefficient vector

$$\mathbf{q}_{101}^{T}=[\mathbf{q}_{101,u}^{T}\quad \mathbf{q}_{101,v}^{T}\quad \mathbf{q}_{101,w}^{T}]. \tag{83}$$

Accordingly, we define the matrix of basis functions as

$$\mathbf{\Phi}_{10,1}^{T}(x_1,x_2)=\mathbf{H}_{10,1}(x_2)\cdot\mathbf{p}_{101,R}^{T}(x_1), \tag{84}$$

where the matrix

$$\mathbf{H}_{10,1}(x_2)=\begin{bmatrix}\dfrac{1}{2} & 0 & 0\\ 0 & \dfrac{1}{2} & 0\\ 0 & 0 & \dfrac{1}{2}\end{bmatrix}, \tag{85}$$

and the matrix of higher order partial derivatives of the basis functions as

$$\mathbf{\Phi}_{10,1}^{(k_1,k_2)T}(x_1,x_2)=\mathbf{H}_{10,1}^{(k_2)}(x_2)\cdot\mathbf{p}_{101,R}^{(k_1)T}(x_1), \tag{86}$$

where $k_1$ and $k_2$ are nonnegative integers, and the function matrixes $\mathbf{p}_{101,R}^{(k_1)T}(x_1)$ and $\mathbf{H}_{10,1}^{(k_2)}(x_2)$ are similar to that of Eqs. (61) and (62).

Therefore, we can define the matrixes of basis functions as

$$\boldsymbol{\Phi}_{1,1}^{\mathrm{T}}(x_1,x_2) = [\boldsymbol{\Phi}_{10,1}^{\mathrm{T}}(x_1,x_2) \quad \boldsymbol{\Phi}_{11,1}^{\mathrm{T}}(x_1,x_2) \quad \cdots \quad \boldsymbol{\Phi}_{1n,1}^{\mathrm{T}}(x_1,x_2) \quad \cdots], \tag{87}$$

$$\boldsymbol{\Phi}_{1,2}^{\mathrm{T}}(x_1,x_2) = [\boldsymbol{\Phi}_{11,2}^{\mathrm{T}}(x_1,x_2) \quad \boldsymbol{\Phi}_{12,2}^{\mathrm{T}}(x_1,x_2) \quad \cdots \quad \boldsymbol{\Phi}_{1n,2}^{\mathrm{T}}(x_1,x_2) \quad \cdots], \tag{88}$$

$$\boldsymbol{\Phi}_{1}^{\mathrm{T}}(x_1,x_2) = [\boldsymbol{\Phi}_{1,1}^{\mathrm{T}}(x_1,x_2) \quad \boldsymbol{\Phi}_{1,2}^{\mathrm{T}}(x_1,x_2)], \tag{89}$$

and define the corresponding matrixes of higher order partial derivatives of the basis functions as

$$\boldsymbol{\Phi}_{1,1}^{(k_1,k_2)\mathrm{T}}(x_1,x_2) = [\boldsymbol{\Phi}_{10,1}^{(k_1,k_2)\mathrm{T}}(x_1,x_2) \quad \boldsymbol{\Phi}_{11,1}^{(k_1,k_2)\mathrm{T}}(x_1,x_2) \quad \cdots \quad \boldsymbol{\Phi}_{1n,1}^{(k_1,k_2)\mathrm{T}}(x_1,x_2) \quad \cdots], \tag{90}$$

$$\boldsymbol{\Phi}_{1,2}^{(k_1,k_2)\mathrm{T}}(x_1,x_2) = [\boldsymbol{\Phi}_{11,2}^{(k_1,k_2)\mathrm{T}}(x_1,x_2) \quad \boldsymbol{\Phi}_{12,2}^{(k_1,k_2)\mathrm{T}}(x_1,x_2) \quad \cdots \quad \boldsymbol{\Phi}_{1n,2}^{(k_1,k_2)\mathrm{T}}(x_1,x_2) \quad \cdots], \tag{91}$$

$$\boldsymbol{\Phi}_{1}^{(k_1,k_2)\mathrm{T}}(x_1,x_2) = [\boldsymbol{\Phi}_{1,1}^{(k_1,k_2)\mathrm{T}}(x_1,x_2) \quad \boldsymbol{\Phi}_{1,2}^{(k_1,k_2)\mathrm{T}}(x_1,x_2)], \tag{92}$$

and define the vectors of boundary Fourier coefficient as

$$\mathbf{q}_{1,1}^{\mathrm{T}} = [\mathbf{q}_{101}^{\mathrm{T}} \quad \mathbf{q}_{111}^{\mathrm{T}} \quad \cdots \quad \mathbf{q}_{1n1}^{\mathrm{T}} \quad \cdots], \tag{93}$$

$$\mathbf{q}_{1,2}^{\mathrm{T}} = [\mathbf{q}_{112}^{\mathrm{T}} \quad \mathbf{q}_{122}^{\mathrm{T}} \quad \cdots \quad \mathbf{q}_{1n2}^{\mathrm{T}} \quad \cdots], \tag{94}$$

$$\mathbf{q}_{1}^{\mathrm{T}} = [\mathbf{q}_{1,1}^{\mathrm{T}} \quad \mathbf{q}_{1,2}^{\mathrm{T}}]. \tag{95}$$

Then the boundary functions $\varphi_{1u}(x_1,x_2)$, $\varphi_{1v}(x_1,x_2)$ and $\varphi_{1w}(x_1,x_2)$ are combined to be expressed as

$$\begin{bmatrix} \varphi_{1u}(x_1,x_2) \\ \varphi_{1v}(x_1,x_2) \\ \varphi_{1w}(x_1,x_2) \end{bmatrix} = \boldsymbol{\Phi}_1^{\mathrm{T}}(x_1,x_2) \cdot \mathbf{q}_1. \tag{96}$$

*3.4. Expressions of boundary functions expanded along the $x_1$-direction*

Suppose that Eq. (8) has the following homogeneous solution

$$\left.\begin{array}{l} \varphi_{2m,u}(x_1,x_2) = \zeta_{1m,u}(x_2)\cos(\alpha_m x_1) + \zeta_{2m,u}(x_2)\sin(\alpha_m x_1) \\ \varphi_{2m,v}(x_1,x_2) = \zeta_{1m,v}(x_2)\cos(\alpha_m x_1) + \zeta_{2m,v}(x_2)\sin(\alpha_m x_1) \\ \varphi_{2m,w}(x_1,x_2) = \zeta_{1m,w}(x_2)\cos(\alpha_m x_1) + \zeta_{2m,w}(x_2)\sin(\alpha_m x_1) \end{array}\right\}, \tag{97}$$

where $m$ is a nonnegative integer, the functions $\zeta_{1m,u}(x_2)$, $\zeta_{2m,u}(x_2)$, $\zeta_{1m,v}(x_2)$, $\zeta_{2m,v}(x_2)$, $\zeta_{1m,w}(x_2)$ and $\zeta_{2m,w}(x_2)$ are undetermined one-dimensional functions, and $\alpha_m = m\pi/a$.

With the same procedure as in section 3.3, we can also obtain the expressions of matrixes of basis functions, i.e., $\boldsymbol{\Phi}_{2m,1}^{\mathrm{T}}(x_1,x_2)$, $\boldsymbol{\Phi}_{2m,2}^{\mathrm{T}}(x_1,x_2)$ and $\boldsymbol{\Phi}_{20,1}^{\mathrm{T}}(x_1,x_2)$.

However, we take advantage of the intrinsic symmetry between the homogeneous solutions (16) and (97), and determine in a direct manner the matrixes of basis functions, $\boldsymbol{\Phi}_{2m,1}^{\mathrm{T}}(x_1,x_2)$, $\boldsymbol{\Phi}_{2m,2}^{\mathrm{T}}(x_1,x_2)$, $\boldsymbol{\Phi}_{20,1}^{\mathrm{T}}(x_1,x_2)$, and their higher order partial derivatives:

1. Based on the expressions of matrixes of basis functions $\boldsymbol{\Phi}_{1n,1}^{\mathrm{T}}(x_1,x_2)$, $\boldsymbol{\Phi}_{1n,2}^{\mathrm{T}}(x_1,x_2)$ and $\boldsymbol{\Phi}_{10,1}^{\mathrm{T}}(x_1,x_2)$, interchange the parameters $n$ with $m$, $b$ with $a$, $x_1$ with $x_2$, and interchange the first row with the second row of the matrix, then interchange the first and second columns with the third and fourth columns of the matrix. Then we obtain the expressions of matrixes of basis functions $\boldsymbol{\Phi}_{2m,1}^{\mathrm{T}}(x_1,x_2)$, $\boldsymbol{\Phi}_{2m,2}^{\mathrm{T}}(x_1,x_2)$ and $\boldsymbol{\Phi}_{20,1}^{\mathrm{T}}(x_1,x_2)$.

2. For nonnegative integers $k_1$ and $k_2$, based on the expressions of matrixes of basis functions $\boldsymbol{\Phi}_{1n,1}^{(k_1,k_2)\mathrm{T}}(x_1,x_2)$, $\boldsymbol{\Phi}_{1n,2}^{(k_1,k_2)\mathrm{T}}(x_1,x_2)$ and $\boldsymbol{\Phi}_{10,1}^{(k_1,k_2)\mathrm{T}}(x_1,x_2)$, we interchange the parameters $n$ with $m$, $b$ with $a$, $x_1$ with $x_2$, $k_1$ with $k_2$, and interchange the first row with the second row of the corresponding matrix, then interchange the first and second columns with the third and fourth columns, and we obtain the expressions of matrixes of basis functions $\mathbf{\Phi}_{2m,1}^{(k_1,k_2)\mathrm{T}}(x_1,x_2)$, $\mathbf{\Phi}_{2m,2}^{(k_1,k_2)\mathrm{T}}(x_1,x_2)$ and $\mathbf{\Phi}_{20,1}^{(k_1,k_2)\mathrm{T}}(x_1,x_2)$.

Therefore, we define the matrixes of basis functions as

$$\mathbf{\Phi}_{2,1}^{\mathrm{T}}(x_1,x_2) = [\mathbf{\Phi}_{20,1}^{\mathrm{T}}(x_1,x_2) \quad \mathbf{\Phi}_{21,1}^{\mathrm{T}}(x_1,x_2) \quad \cdots \quad \mathbf{\Phi}_{2m,1}^{\mathrm{T}}(x_1,x_2) \quad \cdots], \tag{98}$$

$$\mathbf{\Phi}_{2,2}^{\mathrm{T}}(x_1,x_2) = [\mathbf{\Phi}_{21,2}^{\mathrm{T}}(x_1,x_2) \quad \mathbf{\Phi}_{22,2}^{\mathrm{T}}(x_1,x_2) \quad \cdots \quad \mathbf{\Phi}_{2m,2}^{\mathrm{T}}(x_1,x_2) \quad \cdots], \tag{99}$$

$$\mathbf{\Phi}_{2}^{\mathrm{T}}(x_1,x_2) = [\mathbf{\Phi}_{2,1}^{\mathrm{T}}(x_1,x_2) \quad \mathbf{\Phi}_{2,2}^{\mathrm{T}}(x_1,x_2)], \tag{100}$$

and define the corresponding matrixes of higher order partial derivatives of the basis functions as

$$\mathbf{\Phi}_{2,1}^{(k_1,k_2)\mathrm{T}}(x_1,x_2) = [\mathbf{\Phi}_{20,1}^{(k_1,k_2)\mathrm{T}}(x_1,x_2) \quad \mathbf{\Phi}_{21,1}^{(k_1,k_2)\mathrm{T}}(x_1,x_2) \quad \cdots \quad \mathbf{\Phi}_{2m,1}^{(k_1,k_2)\mathrm{T}}(x_1,x_2) \quad \cdots], \tag{101}$$

$$\mathbf{\Phi}_{2,2}^{(k_1,k_2)\mathrm{T}}(x_1,x_2) = [\mathbf{\Phi}_{21,2}^{(k_1,k_2)\mathrm{T}}(x_1,x_2) \quad \mathbf{\Phi}_{22,2}^{(k_1,k_2)\mathrm{T}}(x_1,x_2) \quad \cdots \quad \mathbf{\Phi}_{2m,2}^{(k_1,k_2)\mathrm{T}}(x_1,x_2) \quad \cdots], \tag{102}$$

$$\mathbf{\Phi}_{2}^{(k_1,k_2)\mathrm{T}}(x_1,x_2) = [\mathbf{\Phi}_{2,1}^{(k_1,k_2)\mathrm{T}}(x_1,x_2) \quad \mathbf{\Phi}_{2,2}^{(k_1,k_2)\mathrm{T}}(x_1,x_2)], \tag{103}$$

and define the vectors of boundary Fourier coefficient as

$$\mathbf{q}_{2m1}^{\mathrm{T}} = [\mathbf{q}_{2m2,u}^{\mathrm{T}} \quad \mathbf{q}_{2m1,v}^{\mathrm{T}} \quad \mathbf{q}_{2m1,w}^{\mathrm{T}}], \tag{104}$$

$$\mathbf{q}_{2m2}^{\mathrm{T}} = [\mathbf{q}_{2m1,u}^{\mathrm{T}} \quad \mathbf{q}_{2m2,v}^{\mathrm{T}} \quad \mathbf{q}_{2m2,w}^{\mathrm{T}}], \tag{105}$$

$$\mathbf{q}_{201}^{\mathrm{T}} = [\mathbf{q}_{201,u}^{\mathrm{T}} \quad \mathbf{q}_{201,v}^{\mathrm{T}} \quad \mathbf{q}_{201,w}^{\mathrm{T}}], \tag{106}$$

and

$$\mathbf{q}_{2,1}^{\mathrm{T}} = [\mathbf{q}_{201}^{\mathrm{T}} \quad \mathbf{q}_{211}^{\mathrm{T}} \quad \cdots \quad \mathbf{q}_{2m1}^{\mathrm{T}} \quad \cdots], \tag{107}$$

$$\mathbf{q}_{2,2}^{\mathrm{T}} = [\mathbf{q}_{212}^{\mathrm{T}} \quad \mathbf{q}_{222}^{\mathrm{T}} \quad \cdots \quad \mathbf{q}_{2m2}^{\mathrm{T}} \quad \cdots], \tag{108}$$

$$\mathbf{q}_{2}^{\mathrm{T}} = [\mathbf{q}_{2,1}^{\mathrm{T}} \quad \mathbf{q}_{2,2}^{\mathrm{T}}], \tag{109}$$

then the boundary functions $\varphi_{2u}(x_1,x_2)$, $\varphi_{2v}(x_1,x_2)$ and $\varphi_{2w}(x_1,x_2)$ are combined to be

$$\begin{bmatrix} \varphi_{2u}(x_1,x_2) \\ \varphi_{2v}(x_1,x_2) \\ \varphi_{2w}(x_1,x_2) \end{bmatrix} = \mathbf{\Phi}_{2}^{\mathrm{T}}(x_1,x_2) \cdot \mathbf{q}_{2}. \tag{110}$$

## 3.5. Expressions of internal functions and corner functions

According to the Fourier series multiscale method of the second order linear differential equation with constant coefficients [37], the internal functions of the modal functions $\varphi_u$, $\varphi_v$ and $\varphi_w$ are all full-range Fourier series over the domain $[-a,a]\times[-b,b]$. Without loss of generality, we suppose that

$$\varphi_{0u}(x_1,x_2) = \sum_{m=0}^{\infty}\sum_{n=0}^{\infty} \lambda_{mn}[V_{1mn,u}\cos(\alpha_m x_1)\cos(\beta_n x_2) + V_{2mn,u}\sin(\alpha_m x_1)\cos(\beta_n x_2)$$
$$+ V_{3mn,u}\cos(\alpha_m x_1)\sin(\beta_n x_2) + V_{4mn,u}\sin(\alpha_m x_1)\sin(\beta_n x_2)], \tag{111}$$

$$\varphi_{0v}(x_1,x_2) = \sum_{m=0}^{\infty}\sum_{n=0}^{\infty} \lambda_{mn}[V_{1mn,v}\cos(\alpha_m x_1)\cos(\beta_n x_2) + V_{2mn,v}\sin(\alpha_m x_1)\cos(\beta_n x_2)$$
$$+ V_{3mn,v}\cos(\alpha_m x_1)\sin(\beta_n x_2) + V_{4mn,v}\sin(\alpha_m x_1)\sin(\beta_n x_2)], \tag{112}$$

$$\varphi_{0w}(x_1,x_2) = \sum_{m=0}^{\infty}\sum_{n=0}^{\infty} \lambda_{mn}[V_{1mn,w}\cos(\alpha_m x_1)\cos(\beta_n x_2) + V_{2mn,w}\sin(\alpha_m x_1)\cos(\beta_n x_2)$$

$$+V_{3mn,w}\cos(\alpha_m x_1)\sin(\beta_n x_2)+V_{4mn,w}\sin(\alpha_m x_1)\sin(\beta_n x_2)], \tag{113}$$

where $\alpha_m = m\pi/a$, $\beta_n = n\pi/b$,

$$\lambda_{mn} = \begin{cases} 1/4 & m=0, n=0 \\ 1/2 & m>0, n=0 \\ 1/2 & m=0, n>0 \\ 1 & m>0, n>0 \end{cases},$$

$V_{imn,u}$, $V_{imn,v}$ and $V_{imn,w}$, $i=1,2,3,4$, are the Fourier coefficients of functions $\varphi_{0u}(x_1,x_2)$, $\varphi_{0v}(x_1,x_2)$ and $\varphi_{0w}(x_1,x_2)$ respectively.

We rewrite Eqs. (111)-(113) as

$$\varphi_{0u}(x_1,x_2) = \mathbf{\Phi}_0^{\mathrm{T}}(x_1,x_2)\cdot \mathbf{q}_{0,u}, \tag{114}$$

$$\varphi_{0v}(x_1,x_2) = \mathbf{\Phi}_0^{\mathrm{T}}(x_1,x_2)\cdot \mathbf{q}_{0,v}, \tag{115}$$

$$\varphi_{0w}(x_1,x_2) = \mathbf{\Phi}_0^{\mathrm{T}}(x_1,x_2)\cdot \mathbf{q}_{0,w}, \tag{116}$$

where the vector of basis function $\mathbf{\Phi}_0^{\mathrm{T}}$ and the vectors of Fourier coefficients $\mathbf{q}_{0,u}$, $\mathbf{q}_{0,v}$ and $\mathbf{q}_{0,w}$ are defined as in Section 3.2, [36].

In addition, the corner functions of the modal functions $\varphi_u$, $\varphi_v$ and $\varphi_w$ can be expressed as

$$\varphi_{3u}(x_1,x_2) = \mathbf{\Phi}_3^{\mathrm{T}}(x_1,x_2)\cdot \mathbf{q}_{3,u}, \tag{117}$$

$$\varphi_{3v}(x_1,x_2) = \mathbf{\Phi}_3^{\mathrm{T}}(x_1,x_2)\cdot \mathbf{q}_{3,v}, \tag{118}$$

$$\varphi_{3w}(x_1,x_2) = \mathbf{\Phi}_3^{\mathrm{T}}(x_1,x_2)\cdot \mathbf{q}_{3,w}, \tag{119}$$

where the vector of basis function $\mathbf{\Phi}_3^{\mathrm{T}}$ and the undetermined coefficient vectors $\mathbf{q}_{3,u}$, $\mathbf{q}_{3,v}$ and $\mathbf{q}_{3,w}$ are also defined as in Section 3.2, [36].

It shall be noted that, Eq (8) is a system of second order differential equations. Therefore, the vector of basis function $\mathbf{\Phi}_3^{\mathrm{T}}$ and the undetermined coefficient vectors $\mathbf{q}_{3,u}$, $\mathbf{q}_{3,v}$ and $\mathbf{q}_{3,w}$ contain only one element each. Then we have

$$\mathbf{\Phi}_3^{\mathrm{T}} = [\varphi_{3,00}(x_1,x_2)], \quad \mathbf{q}_{3,u}^{\mathrm{T}} = [q_{3u,00}], \quad \mathbf{q}_{3,v}^{\mathrm{T}} = [q_{3v,00}], \quad \mathbf{q}_{3,w}^{\mathrm{T}} = [q_{3w,00}], \tag{120}$$

where the function

$$\varphi_{3,00}(x_1,x_2) = \frac{x_1 x_2}{4ab}, \tag{121}$$

$q_{3u,00}$, $q_{3v,00}$, $q_{3w,00}$ are undetermined constants.

Substituting Eqs. (114)-(121) into the following equation

$$\mathbf{L}_\varphi \begin{bmatrix} \varphi_{0u}+\varphi_{3u} \\ \varphi_{0v}+\varphi_{3v} \\ \varphi_{0w}+\varphi_{3w} \end{bmatrix} = 0, \tag{122}$$

we then obtain

$$\mathbf{L}_\varphi \begin{bmatrix} \mathbf{\Phi}_0^{\mathrm{T}} & 0 & 0 \\ 0 & \mathbf{\Phi}_0^{\mathrm{T}} & 0 \\ 0 & 0 & \mathbf{\Phi}_0^{\mathrm{T}} \end{bmatrix} \begin{bmatrix} \mathbf{q}_{0,u} \\ \mathbf{q}_{0,v} \\ \mathbf{q}_{0,w} \end{bmatrix} = -\mathbf{L}_\varphi \begin{bmatrix} \mathbf{\Phi}_3^{\mathrm{T}} & 0 & 0 \\ 0 & \mathbf{\Phi}_3^{\mathrm{T}} & 0 \\ 0 & 0 & \mathbf{\Phi}_3^{\mathrm{T}} \end{bmatrix} \begin{bmatrix} \mathbf{q}_{3,u} \\ \mathbf{q}_{3,v} \\ \mathbf{q}_{3,w} \end{bmatrix}. \tag{123}$$

In the forgoing equation, we expand the functions on the right side in full-range Fourier series over the domain $[-a,a] \times [-b,b]$, and compare the Fourier coefficients on both sides successively. Meanwhile, we truncate all the series expansions to $m = M$ and $n = N$ for the sake of numerical implementation, and then establish the relations between the undetermined constant vectors $\mathbf{q}_{0,u}^T$, $\mathbf{q}_{0,v}^T$, $\mathbf{q}_{0,w}^T$ and $\mathbf{q}_{3,u}^T$, $\mathbf{q}_{3,v}^T$, $\mathbf{q}_{3,w}^T$ as follows:

a. if $m = 0$ and $n = 0$, then

$$\left.\begin{aligned}
(\rho\omega^2 - \mu k^2)V_{100,u} &= -(\lambda+\mu)\frac{1}{ab}q_{3v,00} \\
(\rho\omega^2 - \mu k^2)V_{100,v} &= -(\lambda+\mu)\frac{1}{ab}q_{3u,00} \\
[(\rho\omega^2 - (\lambda+2\mu)k^2]V_{100,w} &= 0
\end{aligned}\right\}, \quad (124)$$

b. if $m > 0$ and $n = 0$, then

$$\left.\begin{aligned}
[-(\lambda+2\mu)\alpha_m^2 + \rho\omega^2 - \mu k^2]V_{1m0,u} + [(\lambda+\mu)k\alpha_m]V_{2m0,w} &= 0 \\
[-(\lambda+2\mu)\alpha_m^2 + \rho\omega^2 - \mu k^2]V_{2m0,u} + [-(\lambda+\mu)k\alpha_m]V_{1m0,w} &= 0 \\
[-\mu\alpha_m^2 + \rho\omega^2 - \mu k^2]V_{1m0,v} &= 0 \\
[-\mu\alpha_m^2 + \rho\omega^2 - \mu k^2]V_{2m0,v} &= (\lambda+\mu)k\frac{(-1)^m}{bm\pi}q_{3w,00} \\
[-(\lambda+\mu)k\alpha_m]V_{2m0,u} + [-\mu\alpha_m^2 + \rho\omega^2 - (\lambda+2\mu)k^2]V_{1m0,w} &= 0 \\
[(\lambda+\mu)k\alpha_m]V_{1m0,u} + [-\mu\alpha_m^2 + \rho\omega^2 - (\lambda+2\mu)k^2]V_{2m0,w} & \\
= (\lambda+\mu)k\frac{(-1)^{m+1}}{bm\pi}q_{3v,00} &
\end{aligned}\right\}, \quad (125)$$

c. if $m = 0$ and $n > 0$, then

$$\left.\begin{aligned}
[-\mu\beta_n^2 + \rho\omega^2 - \mu k^2]V_{10n,u} &= 0 \\
[-\mu\beta_n^2 + \rho\omega^2 - \mu k^2]V_{30n,u} &= (\lambda+\mu)k\frac{(-1)^n}{an\pi}q_{3w,00} \\
[-(\lambda+2\mu)\beta_n^2 + \rho\omega^2 - \mu k^2]V_{10n,v} + [(\lambda+\mu)k\beta_n]V_{30n,w} &= 0 \\
[-(\lambda+2\mu)\beta_n^2 + \rho\omega^2 - \mu k^2]V_{30n,v} + [-(\lambda+\mu)k\beta_n]V_{10n,w} &= 0 \\
[-(\lambda+\mu)k\beta_n]V_{30n,v} + [-\mu\beta_n^2 + \rho\omega^2 - (\lambda+2\mu)k^2]V_{10n,w} &= 0 \\
[(\lambda+\mu)k\beta_n]V_{10n,v} + [-\mu\beta_n^2 + \rho\omega^2 - (\lambda+2\mu)k^2]V_{30n,w} & \\
= (\lambda+\mu)k\frac{(-1)^{n+1}}{an\pi}q_{3u,00} &
\end{aligned}\right\}, \quad (126)$$

d. if $m > 0$ and $n > 0$, then

$$\left.\begin{aligned}&[-(\lambda+2\mu)\alpha_m^2-\mu\beta_n^2+\rho\omega^2-\mu k^2]V_{1mn,u}+[(\lambda+\mu)\alpha_m\beta_n]V_{4mn,v}\\&+[(\lambda+\mu)k\alpha_m]V_{2mn,w}=0\\&[(\lambda+\mu)\alpha_m\beta_n]V_{1mn,u}+[-\mu\alpha_m^2-(\lambda+2\mu)\beta_n^2+\rho\omega^2-\mu k^2]V_{4mn,v}\\&+[-(\lambda+\mu)k\beta_n]V_{2mn,w}=(\rho\omega^2-\mu k^2)\frac{(-1)^{m+n+1}}{mn\pi^2}q_{3v,00}\\&[(\lambda+\mu)k\alpha_m]V_{1mn,u}+[-(\lambda+\mu)k\beta_n]V_{4mn,v}\\&+[-\mu\alpha_m^2-\mu\beta_n^2+\rho\omega^2-(\lambda+2\mu)k^2]V_{2mn,w}=0\end{aligned}\right\}, \quad (127)$$

$$\left.\begin{aligned}&[-(\lambda+2\mu)\alpha_m^2-\mu\beta_n^2+\rho\omega^2-\mu k^2]V_{2mn,u}+[-(\lambda+\mu)\alpha_m\beta_n]V_{3mn,v}\\&+[-(\lambda+\mu)k\alpha_m]V_{1mn,w}=0\\&[-(\lambda+\mu)\alpha_m\beta_n]V_{2mn,u}+[-\mu\alpha_m^2-(\lambda+2\mu)\beta_n^2+\rho\omega^2-\mu k^2]V_{3mn,v}\\&+[-(\lambda+\mu)k\beta_n]V_{1mn,w}=0\\&[-(\lambda+\mu)k\alpha_m]V_{2mn,u}+[-(\lambda+\mu)k\beta_n]V_{3mn,v}\\&+[-\mu\alpha_m^2-\mu\beta_n^2+\rho\omega^2-(\lambda+2\mu)k^2]V_{1mn,w}=0\end{aligned}\right\}, \quad (128)$$

$$\left.\begin{aligned}&[-(\lambda+2\mu)\alpha_m^2-\mu\beta_n^2+\rho\omega^2-\mu k^2]V_{3mn,u}+[-(\lambda+\mu)\alpha_m\beta_n]V_{2mn,v}\\&+[(\lambda+\mu)k\alpha_m]V_{4mn,w}=0\\&[-(\lambda+\mu)\alpha_m\beta_n]V_{3mn,u}+[-\mu\alpha_m^2-(\lambda+2\mu)\beta_n^2+\rho\omega^2-\mu k^2]V_{2mn,v}\\&+[(\lambda+\mu)k\beta_n]V_{4mn,w}=0\\&[(\lambda+\mu)k\alpha_m]V_{3mn,u}+[(\lambda+\mu)k\beta_n]V_{2mn,v}\\&+[-\mu\alpha_m^2-\mu\beta_n^2+\rho\omega^2-(\lambda+2\mu)k^2]V_{4mn,w}\\&=[\rho\omega^2-(\lambda+2\mu)k^2]\frac{(-1)^{m+n+1}}{mn\pi^2}q_{3w,00}\end{aligned}\right\}, \quad (129)$$

and

$$\left.\begin{aligned}&[-(\lambda+2\mu)\alpha_m^2-\mu\beta_n^2+\rho\omega^2-\mu k^2]V_{4mn,u}+[(\lambda+\mu)\alpha_m\beta_n]V_{1mn,v}\\&+[-(\lambda+\mu)k\alpha_m]V_{3mn,w}=(\rho\omega^2-\mu k^2)\frac{(-1)^{m+n+1}}{mn\pi^2}q_{3u,00}\\&[(\lambda+\mu)\alpha_m\beta_n]V_{4mn,u}+[-\mu\alpha_m^2-(\lambda+2\mu)\beta_n^2+\rho\omega^2-\mu k^2]V_{1mn,v}\\&+[(\lambda+\mu)k\beta_n]V_{3mn,w}=0\\&[-(\lambda+\mu)k\alpha_m]V_{4mn,u}+[(\lambda+\mu)k\beta_n]V_{1mn,v}\\&+[-\mu\alpha_m^2-\mu\beta_n^2+\rho\omega^2-(\lambda+2\mu)k^2]V_{3mn,w}=0\end{aligned}\right\}. \quad (130)$$

If we write the undetermined constant vectors as
$$\mathbf{q}_0^T=[\mathbf{q}_{0,u}^T \quad \mathbf{q}_{0,v}^T \quad \mathbf{q}_{0,w}^T], \quad (131)$$
$$\mathbf{q}_3^T=[\mathbf{q}_{3,u}^T \quad \mathbf{q}_{3,v}^T \quad \mathbf{q}_{3,w}^T]=[q_{3u,00} \quad q_{3v,00} \quad q_{3w,00}], \quad (132)$$
then Eqs. (124)-(130) are combined to be
$$\mathbf{q}_0=\mathbf{T}_{03}\mathbf{q}_3, \quad (133)$$
where $\mathbf{T}_{03}$ is a transformation matrix.

Accordingly, the internal functions and corner functions are expressed as

$$\begin{bmatrix} \varphi_{0u}(x_1,x_2)+\varphi_{3u}(x_1,x_2) \\ \varphi_{0v}(x_1,x_2)+\varphi_{3v}(x_1,x_2) \\ \varphi_{0w}(x_1,x_2)+\varphi_{3w}(x_1,x_2) \end{bmatrix} = \mathbf{\Phi}_{03,R}^{\mathrm{T}}(x_1,x_2)\cdot\mathbf{q}_3, \tag{134}$$

where the matrix

$$\mathbf{\Phi}_{03,R}^{\mathrm{T}}(x_1,x_2) = \begin{bmatrix} \mathbf{\Phi}_0^{\mathrm{T}} & \mathbf{0} & \mathbf{0} \\ \mathbf{0} & \mathbf{\Phi}_0^{\mathrm{T}} & \mathbf{0} \\ \mathbf{0} & \mathbf{0} & \mathbf{\Phi}_0^{\mathrm{T}} \end{bmatrix} \cdot \mathbf{T}_{03} + \begin{bmatrix} \mathbf{\Phi}_3^{\mathrm{T}} & \mathbf{0} & \mathbf{0} \\ \mathbf{0} & \mathbf{\Phi}_3^{\mathrm{T}} & \mathbf{0} \\ \mathbf{0} & \mathbf{0} & \mathbf{\Phi}_3^{\mathrm{T}} \end{bmatrix}. \tag{135}$$

*3.6. Expression of the Fourier series multiscale solution*

If we write the matrix of basis functions

$$\mathbf{\Phi}^{\mathrm{T}}(x_1,x_2) = [\mathbf{\Phi}_{03,R}^{\mathrm{T}}(x_1,x_2)\quad \mathbf{\Phi}_1^{\mathrm{T}}(x_1,x_2)\quad \mathbf{\Phi}_2^{\mathrm{T}}(x_1,x_2)], \tag{136}$$

and the undetermined constant vector

$$\mathbf{q}^{\mathrm{T}} = [\mathbf{q}_3^{\mathrm{T}}\quad \mathbf{q}_1^{\mathrm{T}}\quad \mathbf{q}_2^{\mathrm{T}}], \tag{137}$$

then the full-range composite Fourier series of modal functions $\varphi_u$, $\varphi_v$ and $\varphi_w$ over the domain $[-a,a]\times[-b,b]$ can be expressed as

$$\begin{bmatrix} \varphi_u(x_1,x_2) \\ \varphi_v(x_1,x_2) \\ \varphi_w(x_1,x_2) \end{bmatrix} = \mathbf{\Phi}^{\mathrm{T}}(x_1,x_2)\cdot\mathbf{q}. \tag{138}$$

We write the matrix $\mathbf{\Phi}^{\mathrm{T}}(x_1,x_2)$ in the following block form

$$\mathbf{\Phi}^{\mathrm{T}}(x_1,x_2) = \begin{bmatrix} \mathbf{\Phi}_u^{\mathrm{T}}(x_1,x_2) \\ \mathbf{\Phi}_v^{\mathrm{T}}(x_1,x_2) \\ \mathbf{\Phi}_w^{\mathrm{T}}(x_1,x_2) \end{bmatrix}, \tag{139}$$

and by combining Eq. (7), we express the Fourier series multiscale solution for wave propagation in a rectangular beam as

$$\left.\begin{aligned} u &= [\mathbf{\Phi}_u^{\mathrm{T}}(x_1,x_2)\cdot\mathbf{q}]\cos(kz-\omega t) \\ v &= [\mathbf{\Phi}_v^{\mathrm{T}}(x_1,x_2)\cdot\mathbf{q}]\cos(kz-\omega t) \\ w &= [\mathbf{\Phi}_w^{\mathrm{T}}(x_1,x_2)]\cdot\mathbf{q}]\sin(kz-\omega t) \end{aligned}\right\}. \tag{140}$$

*3.7. Expressions of stress resultants*

Substituting the above Fourier series multiscale solution into the constitutive equation (4), we obtain the normal stresses $\sigma_{x_1}$, $\sigma_{x_2}$, $\sigma_z$ and shear stresses $\tau_{x_1 x_2}$, $\tau_{x_2 z}$, $\tau_{x_1 z}$ in the beam, and express them as follows

$$\left.\begin{array}{l}\sigma_{x_1} = [\mathbf{\Gamma}_1 \cdot \mathbf{q}]\cos(kz-\omega t) \\ \sigma_{x_2} = [\mathbf{\Gamma}_2 \cdot \mathbf{q}]\cos(kz-\omega t) \\ \sigma_z = [\mathbf{\Gamma}_3 \cdot \mathbf{q}]\cos(kz-\omega t) \\ \tau_{x_1 x_2} = [\mathbf{\Gamma}_4 \cdot \mathbf{q}]\cos(kz-\omega t) \\ \tau_{x_2 z} = [\mathbf{\Gamma}_5 \cdot \mathbf{q}]\sin(kz-\omega t) \\ \tau_{x_1 z} = [\mathbf{\Gamma}_6 \cdot \mathbf{q}]\sin(kz-\omega t)\end{array}\right\}, \tag{141}$$

where the matrix

$$\mathbf{\Gamma} = \begin{bmatrix}\mathbf{\Gamma}_1 \\ \mathbf{\Gamma}_2 \\ \mathbf{\Gamma}_3 \\ \mathbf{\Gamma}_4 \\ \mathbf{\Gamma}_5 \\ \mathbf{\Gamma}_6\end{bmatrix} = \begin{bmatrix}(\lambda+2\mu)\mathbf{\Phi}_u^{(1,0)\mathrm{T}}(x_1,x_2) + \lambda\mathbf{\Phi}_v^{(0,1)\mathrm{T}}(x_1,x_2) + \lambda k\mathbf{\Phi}_w^{\mathrm{T}}(x_1,x_2) \\ \lambda\mathbf{\Phi}_u^{(1,0)\mathrm{T}}(x_1,x_2) + (\lambda+2\mu)\mathbf{\Phi}_v^{(0,1)\mathrm{T}}(x_1,x_2) + \lambda k\mathbf{\Phi}_w^{\mathrm{T}}(x_1,x_2) \\ \lambda\mathbf{\Phi}_u^{(1,0)\mathrm{T}}(x_1,x_2) + \lambda\mathbf{\Phi}_v^{(0,1)\mathrm{T}}(x_1,x_2) + (\lambda+2\mu)k\mathbf{\Phi}_w^{\mathrm{T}}(x_1,x_2) \\ \mu\mathbf{\Phi}_u^{(0,1)\mathrm{T}}(x_1,x_2) + \mu\mathbf{\Phi}_v^{(1,0)\mathrm{T}}(x_1,x_2) \\ -\mu k\mathbf{\Phi}_v^{\mathrm{T}}(x_1,x_2) + \mu\mathbf{\Phi}_w^{(0,1)\mathrm{T}}(x_1,x_2) \\ -\mu k\mathbf{\Phi}_u^{\mathrm{T}}(x_1,x_2) + \mu\mathbf{\Phi}_w^{(1,0)\mathrm{T}}(x_1,x_2)\end{bmatrix}. \tag{142}$$

## 4. Problem solving

When the Fourier series multiscale solution is employed for the analysis of wave propagation in a rectangular beam, we respectively retain $M$ and $N$ terms along the $x_1$- and $x_2$-direction, then the number of variables contained in the undetermined constant vector $\mathbf{q}$ is $12M+12N+15$. If we also choose an appropriate number of collocation points on the boundary, and set each of the truncated series representing the displacement components or the normal and tangential components of stress at each point equal to the specified boundary conditions, we obtain $12M+12N+15$ homogeneous, linear equations for the undetermined constants. When the determinant of this set of equations is equal to zero, we find a transcendental equation relating the wave frequency $\omega$ and the wave number $k$. We usually call this equation the frequency equation. For a fixed wave frequency, the frequency equation determines some permissible wave numbers of the waves propagating in the beam. And with the frequency changing, these wave numbers also change, to form the dispersion curves. In the $\omega-k$ plane, the frequency equation gives uncountable branches of the dispersion curves, where different branches correspond to different deformation patterns of the waveguide, namely the wave modes, and the intersection points of the branches and the $\omega$ axis are the cut-off frequencies.

Similar to the transverse vibration analysis of thin plates [38], symmetric decomposition for the wave modes in the beam of rectangular cross section is performed herein based on symmetry of boundary conditions:

1. As shown in Figure 2.a, the boundary condition for the rectangular beam is three clamped edges and one free edge (CCFC), which is symmetric about the $x_1$-axis. Therefore, we separate the wave modes in the beam into two types, namely, the longitudinal wave ($L$) and the bending wave ($B_{x1}$) about the $x_1$-axis.

2. As shown in Figures 2.b and 2.c, the boundary conditions for the rectangular beam are four free edges (FFFF) or two opposite free edges and two opposite clamped edges (FCFC), which are twofold symmetric about the $x_1$-axis and the $x_2$-axis simultaneously. Therefore, we separate the wave modes in the beam into four types, namely, the torsional wave ($T$), the longitudinal wave ($L$), the bending wave ($B_{x1}$) about the $x_1$-axis and the bending wave ($B_{x2}$) about the $x_2$-axis. As to the beam with square cross section, the first two kinds of waves with FFFF boundary conditions are further separated into two cases, which are symmetric ($T_s$, $L_s$) or asymmetric ($T_a$, $L_a$) with respect to the diagonal $x_1 = x_2$.

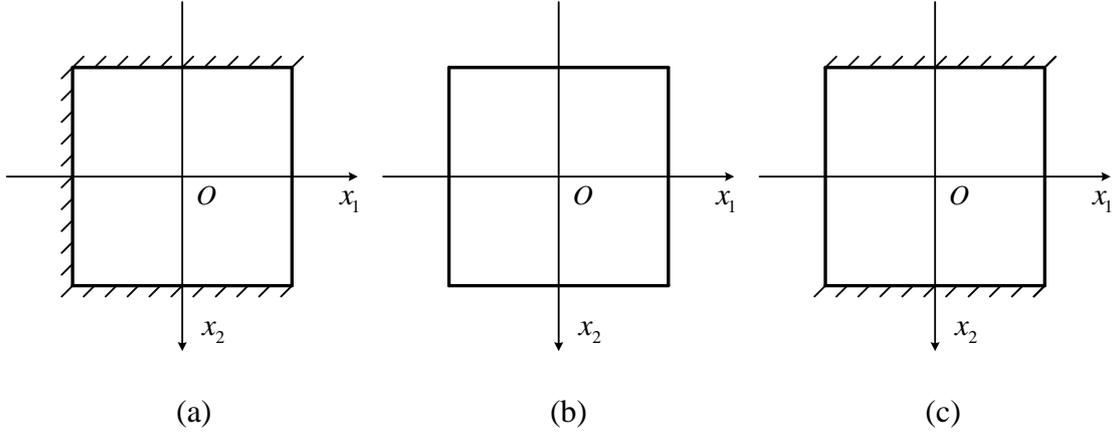

Figure 2: Typical boundary conditions: (a) CCFC, (b) FFFF, (c) FCFC.

For different types of wave modes, we list in Table 2 the symmetric decomposition of displacements in the rectangular beam and the corresponding numbers of undetermined constants involved in the frequency equation.

Table 2: Symmetric decomposition of displacements in a rectangular beam.

| No. | BC | Wave | with respect to the $x_1$-axis | | | with respect to the $x_2$-axis | | | numbers of undetermined constants |
|---|---|---|---|---|---|---|---|---|---|
| | | | $u$ | $v$ | $w$ | $u$ | $v$ | $w$ | |
| 1 | CCFC | $L$ | S | A | S | | | | $6M + 6N + 7$ |
| | | $B_{x1}$ | A | S | A | | | | $6M + 6N + 5$ |
| 2 | FFFF | $L$ | S | A | S | A | S | S | $3M + 3N + 4$ |
| | | $T$ | A | S | A | S | A | A | $3M + 3N + 3$ |
| | | $B_{x1}$ | A | S | A | A | S | S | $3M + 3N + 4$ |
| | | $B_{x2}$ | S | A | S | S | A | A | $3M + 3N + 4$ |
| 3 | FCFC | $L$ | S | A | S | A | S | S | $3M + 3N + 4$ |
| | | $T$ | A | S | A | S | A | A | $3M + 3N + 2$ |
| | | $B_{x1}$ | A | S | A | A | S | S | $3M + 3N + 3$ |
| | | $B_{x2}$ | S | A | S | S | A | A | $3M + 3N + 3$ |

BC: Boundary Condition, S: symmetric, A: asymmetric

# 5. Convergence characteristics

As shown in Table 3, we have specially designed two sets of comparative numerical experiments for the purpose of a detailed analysis of the influences of some key factors, such as computational parameters and boundary conditions, on the convergence characteristics and approximation precision of the Fourier series multiscale solution of the wave propagation in a beam of rectangular cross section. We perform the first set of numerical experiment based on the reference computational scheme, and adjust length-to-width ratios of rectangular cross section of the beam with from 1.0 to 0.67, 0.50, 1.25 and 2.0 successively. We perform the second set of numerical experiments based on the reference computational scheme, and adjust the boundary conditions of wave propagation in a beam of rectangular cross section from four free edges (FFFF) to two opposite free edges and two other opposite clamped edges (FCFC) and three clamped edges and one free edge (CCFC) successively.

It is necessary to point out that the computation is based on the symmetric decomposition of displacements as tabulated in Table 2. The Poisson's ratio is set to be $\upsilon = 0.3$, the non-dimensional wave frequency $\Omega$ and wave number $K$ are defined by $\Omega = \omega a/\pi c_T$ and $K = ka/\pi$, where $c_T = \sqrt{\mu/\rho}$.

Table 3: Computational schemes for wave propagation in a rectangular beam.

| Numerical experiment | No. | Boundary condition | Length-to-width ratio $a/b$ | Derivation technique for discrete equations |
|---|---|---|---|---|
| 1 | a | | 1.0 | |
| | b | | 0.67 | |
| | c | FFFF | 0.50 | the collocation method |
| | d | | 1.25 | |
| | e | | 2.0 | |
| 2 | a | FFFF | | |
| | b | FCFC | 1.0 | the collocation method |
| | c | CCFC | | |

## 5.1. General convergence characteristics

Based on the reference computational schemes as shown in Table 3, the dynamic behavior of the beam is simulated with the truncated composite Fourier series of the Fourier series multiscale solution retaining the first $M$ and $N$ terms, respectively, along the $x_1$- and $x_2$-direction. When $M = N = 4$, 8, 12, 16 and 20, the dispersion curves $\Omega_{M,N}^l(K) - K$ for the first three kinds of waves, altogether five cases, are obtained by using the corresponding frequency equations, where $l$ is the order of the specific traveling wave in the beam. If the computed dispersion curve corresponding to $M = N = 20$ is taken as the reference, the computational error of dispersion curves over the interval of wave number $K \in [0,1]$ is defined as

$$e(\Omega^l_{M,N}) = \max_{K \in \Upsilon^t} \left| \Omega^l_{M,N}(K) - \Omega^l_{20,20}(K) \right|, \tag{143}$$

where the set of sampling points is

$$\Upsilon^t = \{(n_1 - 1)/100, n_1 = 1, 2, \cdots, 101\}. \tag{144}$$

With the increase of the number of truncated terms of the composite Fourier series, the numerical results for the first branches of the dispersion curves and the computational errors of different branches of dispersion curves over the interval of wave number $K \in [0,1]$, as shown in Tables 4 and 5, are displayed as follows:

a. The first branches of dispersion curves of various kinds of wave propagation are uniformly convergent. For the dispersion curves $L^1_s$, $T^1_a$ and $B^1_{x_1}$, when the number of truncated terms equals to 4, the computational errors are already about 4.0E-3 and show a trend of sustained decrease. And when the number of truncated terms equals to 16, the computational errors further decrease by one order of magnitude. For the dispersion curves $L^1_a$ and $T^1_s$, when the number of truncated terms equals to 4, the computational errors are relatively large at about 4.0E-2 and show a trend of sustained and rapid decrease. And when the number of terms equals to 16, the computational errors further decrease to 3.0E-3.

b. The higher order branches of dispersion curves of various kinds of wave propagation are uniformly convergent. For most of the dispersion curves (rather than $L^3_a$ and $L^4_a$), when the number of truncated terms equals to 4, the computational errors are usually between 1.0E-4 and 4.0E-2, and show a trend of sustained decrease. And when the number of truncated terms equals to 16, the computational errors further decrease by 0-2 orders of magnitude. And for the dispersion curves $L^3_a$ and $L^4_a$, when the number of truncated terms equals to 4, the computational errors are relatively large at about 4.0E-2, and show a trend of slow down decrease. And when the number of truncated terms equals to 16, the computational errors are about 1.0E-2.

Table 4: Uniform convergence of first order dispersion curves.

| Wave | K | Ω | | | | | Reference [31] |
|---|---|---|---|---|---|---|---|
| | | M=N=4 | M=N=8 | M=N=12 | M=N=16 | M=N=20 | |
| $L_s$ | 0.3183 | 0.4937 | 0.4937 | 0.4937 | 0.4937 | 0.4937 | 0.4938 |
| | 0.5730 | 0.7668 | 0.7668 | 0.7668 | 0.7668 | 0.7668 | 0.7669 |
| | 0.8276 | 0.9339 | 0.9339 | 0.9339 | 0.9339 | 0.9339 | 0.9340 |
| $L_a$ | 0.3183 | 0.5056 | 0.5261 | 0.5163 | 0.5114 | 0.5084 | 0.6890 |
| | 0.5730 | 0.5958 | 0.5815 | 0.5765 | 0.5740 | 0.5724 | 0.7273 |
| | 0.8276 | 0.6908 | 0.6841 | 0.6818 | 0.6806 | 0.6798 | 0.8490 |
| $T_s$ | 0.3183 | 0.6818 | 0.6721 | 0.6679 | 0.6654 | 0.6637 | 0.6264 |
| | 0.5730 | 0.8058 | 0.7915 | 0.7857 | 0.7822 | 0.7799 | 0.7392 |
| | 0.8276 | 0.9818 | 0.9664 | 0.9604 | 0.9570 | 0.9546 | 0.9177 |
| $T_a$ | 0.3183 | 0.2926 | 0.2923 | 0.2922 | 0.2922 | 0.2922 | 0.2923 |
| | 0.5730 | 0.5263 | 0.5257 | 0.5256 | 0.5255 | 0.5255 | 0.5255 |
| | 0.8276 | 0.7594 | 0.7585 | 0.7582 | 0.7581 | 0.7581 | 0.7581 |
| $B_{x1}$ | 0.3183 | 0.2047 | 0.2038 | 0.2035 | 0.2033 | 0.2032 | 0.2010 |
| | 0.5730 | 0.4575 | 0.4561 | 0.4555 | 0.4552 | 0.4550 | 0.4511 |
| | 0.8276 | 0.7124 | 0.7105 | 0.7097 | 0.7093 | 0.7090 | 0.7032 |

Table 5: Computational error of dispersion curves $e(\Omega^l_{M,N})$.

| Wave | Branch Order $l$ | $M=N=4$ | $M=N=8$ | $M=N=12$ | $M=N=16$ |
|---|---|---|---|---|---|
| $L_s$ | 1 | 0.0001 | 0.0001 | 0.0001 | 0.0000 |
| | 2 | 0.0014 | 0.0002 | 0.0001 | 0.0001 |
| | 3 | 0.0008 | 0.0003 | 0.0002 | 0.0001 |
| | 4 | 0.0013 | 0.0004 | 0.0001 | 0.0001 |
| | 5 | 0.0019 | 0.0003 | 0.0001 | 0.0001 |
| | 6 | 0.0010 | 0.0002 | 0.0001 | 0.0001 |
| | 7 | 0.0141 | 0.0004 | 0.0001 | 0.0001 |
| | 8 | 0.0004 | 0.0004 | 0.0002 | 0.0001 |
| | 9 | 0.0064 | 0.0013 | 0.0004 | 0.0001 |
| $L_a$ | 1 | 0.0397 | 0.0190 | 0.0086 | 0.0033 |
| | 2 | 0.0138 | 0.0047 | 0.0021 | 0.0008 |
| | 3 | 0.0161 | 0.0397 | 0.0273 | 0.0117 |
| | 4 | 0.0392 | 0.0395 | 0.0211 | 0.0127 |
| | 5 | 0.0087 | 0.0022 | 0.0008 | 0.0003 |
| | 6 | 0.0022 | 0.0009 | 0.0004 | 0.0002 |
| | 7 | 0.0396 | 0.0104 | 0.0038 | 0.0014 |
| | 8 | 0.0099 | 0.0027 | 0.0012 | 0.0005 |
| | 9 | 0.0000 | 0.0000 | 0.0000 | 0.0000 |
| $T_s$ | 1 | 0.0272 | 0.0119 | 0.0059 | 0.0024 |
| | 2 | 0.0388 | 0.0325 | 0.0161 | 0.0065 |
| | 3 | 0.0000 | 0.0000 | 0.0000 | 0.0000 |
| | 4 | 0.0105 | 0.0046 | 0.0022 | 0.0009 |
| | 5 | 0.0392 | 0.0393 | 0.0193 | 0.0078 |
| | 6 | 0.0395 | 0.0204 | 0.0096 | 0.0038 |
| | 7 | 0.0395 | 0.0188 | 0.0090 | 0.0036 |
| | 8 | 0.0114 | 0.0071 | 0.0034 | 0.0014 |
| | 9 | 0.0384 | 0.0247 | 0.0119 | 0.0048 |
| | 10 | 0.0396 | 0.0267 | 0.0135 | 0.0054 |
| $T_a$ | 1 | 0.0016 | 0.0005 | 0.0002 | 0.0001 |
| | 2 | 0.0007 | 0.0007 | 0.0003 | 0.0001 |
| | 3 | 0.0034 | 0.0002 | 0.0001 | 0.0001 |
| | 4 | 0.0001 | 0.0001 | 0.0001 | 0.0001 |
| | 5 | 0.0021 | 0.0001 | 0.0001 | 0.0001 |
| | 6 | 0.0006 | 0.0002 | 0.0001 | 0.0001 |
| | 7 | 0.0143 | 0.0003 | 0.0002 | 0.0001 |
| | 8 | 0.0171 | 0.0007 | 0.0001 | 0.0001 |
| | 9 | 0.0090 | 0.0004 | 0.0001 | 0.0001 |
| | 10 | 0.0022 | 0.0002 | 0.0001 | 0.0001 |
| | 11 | 0.0001 | 0.0000 | 0.0000 | 0.0000 |
| $B_{x1}$ | 1 | 0.0041 | 0.0019 | 0.0010 | 0.0004 |
| | 2 | 0.0398 | 0.0218 | 0.0112 | 0.0046 |
| | 3 | 0.0396 | 0.0214 | 0.0110 | 0.0046 |
| | 4 | 0.0030 | 0.0007 | 0.0003 | 0.0001 |
| | 5 | 0.0147 | 0.0072 | 0.0038 | 0.0016 |
| | 6 | 0.0250 | 0.0194 | 0.0121 | 0.0059 |

| | | | | |
|---|---|---|---|---|
| 7  | 0.0390 | 0.0336 | 0.0166 | 0.0067 |
| 8  | 0.0160 | 0.0059 | 0.0028 | 0.0011 |
| 9  | 0.0168 | 0.0082 | 0.0041 | 0.0017 |
| 10 | 0.0083 | 0.0044 | 0.0022 | 0.0009 |
| 11 | 0.0060 | 0.0023 | 0.0011 | 0.0005 |
| 12 | 0.0047 | 0.0027 | 0.0014 | 0.0006 |
| 13 | 0.0251 | 0.0125 | 0.0064 | 0.0027 |
| 14 | 0.0391 | 0.0293 | 0.0145 | 0.0059 |
| 15 | 0.0364 | 0.0169 | 0.0084 | 0.0034 |
| 16 | 0.0138 | 0.0046 | 0.0021 | 0.0008 |
| 17 | 0.0393 | 0.0223 | 0.0131 | 0.0053 |
| 18 | 0.0218 | 0.0114 | 0.0058 | 0.0023 |

Additional remarks about the numerical results in Tables 4 and 5 are as follows:

a. The numerical results of the dispersion curves $L_s^1$, $T_a^1$, $B_{x_1}^1$ are completely consistent with those in [31], while the numerical results of the dispersion curves $L_a^1$ and $T_s^1$ are obviously different from those in [31]. This is related to the selected forms of trial functions in the computational scheme. The trial function we used is the Fourier series multiscale solution of the wave propagation in the beam of rectangular cross section. The trial function is suitable to the rectangular solution domain, satisfies sufficient conditions for two times term-by-term differentiation, and is capable of adjusting adaptively the scales of analysis. However, the trial function used in [31] is the analytical solution of the wave propagation in a beam of circular cross section. This function is actually not suitable to the rectangular solution domain and does not satisfy the sufficient conditions of second order termwise differentiation. Accordingly, the computational scheme used in [31] can only find the first branches of the dispersion curves $L_s^1$, $T_a^1$ and $B_{x_1}^1$ of wave propagation in the beam, and cannot exactly find other branches of dispersion curves, including $L_a^1$, $T_s^1$ and higher order branches of dispersion curves, of wave propagation in the beam.

b. The computational errors of dispersion curves $L_a^3$ and $L_a^4$ are relatively large. This results in difficulties in recognition of $L_a^3$ and $L_a^4$, and even blind area at the position of intersection of these two dispersion curves.

*5.2. Comparative numerical experiments*

1. Numerical experiment 1

We perform the first set of numerical experiment by adjusting length-to-width ratio of wave propagation in a rectangular beam from 1.0 to 0.67, 0.50, 1.25 and 2.0 successively. The dispersion curve $B_{x_1}^1$ is obtained by the Fourier series multiscale solution as an example, and the numerical results are presented in Table 6. A brief discussion is as follows:

a. The dispersion curves of wave propagation are uniformly convergent. When the number of truncated series terms equals to 4, the computational errors are already about 5.0E-3 and show a trend of sustained decrease. When the number of truncated terms equals to 16, the computational errors further decrease by one order of magnitude.

b. The change of length-to-width ratio has little effect on the uniform convergence of the dispersion curves of wave propagation.

Table 6: Uniform convergence of first order dispersion curves $B_{x_1}^1$.

| Length-to-width ratio $a/b$ | $K$ | $\Omega$ | | | | |
|---|---|---|---|---|---|---|
| | | $M=N=4$ | $M=N=8$ | $M=N=12$ | $M=N=16$ | $M=N=20$ |
| 1.0 | 0.1000 | 0.0282 | 0.0280 | 0.0280 | 0.0279 | 0.0279 |
| | 0.3000 | 0.1873 | 0.1865 | 0.1862 | 0.1860 | 0.1859 |
| | 0.5000 | 0.3842 | 0.3829 | 0.3824 | 0.3821 | 0.3819 |
| | 0.7000 | 0.5851 | 0.5834 | 0.5828 | 0.5824 | 0.5821 |
| | 0.9000 | 0.7842 | 0.7822 | 0.7813 | 0.7808 | 0.7805 |
| 0.67 | 0.1000 | 0.0394 | 0.0392 | 0.0392 | 0.0391 | 0.0391 |
| | 0.3000 | 0.2220 | 0.2213 | 0.2211 | 0.2209 | 0.2208 |
| | 0.5000 | 0.4224 | 0.4214 | 0.4210 | 0.4208 | 0.4206 |
| | 0.7000 | 0.6201 | 0.6188 | 0.6182 | 0.6178 | 0.6176 |
| | 0.9000 | 0.8144 | 0.8126 | 0.8118 | 0.8113 | 0.8109 |
| 0.50 | 0.1000 | 0.0486 | 0.0484 | 0.0483 | 0.0483 | 0.0483 |
| | 0.3000 | 0.2413 | 0.2408 | 0.2406 | 0.2405 | 0.2404 |
| | 0.5000 | 0.4402 | 0.4393 | 0.4389 | 0.4387 | 0.4386 |
| | 0.7000 | 0.6345 | 0.6332 | 0.6326 | 0.6322 | 0.6320 |
| | 0.9000 | 0.8253 | 0.8236 | 0.8227 | 0.8221 | 0.8217 |
| 1.25 | 0.1000 | 0.0231 | 0.0230 | 0.0229 | 0.0229 | 0.0229 |
| | 0.3000 | 0.1658 | 0.1649 | 0.1646 | 0.1644 | 0.1643 |
| | 0.5000 | 0.3563 | 0.3549 | 0.3543 | 0.3540 | 0.3538 |
| | 0.7000 | 0.5567 | 0.5549 | 0.5541 | 0.5537 | 0.5534 |
| | 0.9000 | 0.7576 | 0.7554 | 0.7545 | 0.7540 | 0.7536 |
| 2.0 | 0.1000 | 0.0149 | 0.0148 | 0.0148 | 0.0147 | 0.0147 |
| | 0.3000 | 0.1204 | 0.1196 | 0.1193 | 0.1191 | 0.1189 |
| | 0.5000 | 0.2850 | 0.2835 | 0.2829 | 0.2825 | 0.2823 |
| | 0.7000 | 0.4743 | 0.4723 | 0.4714 | 0.4709 | 0.4705 |
| | 0.9000 | 0.6728 | 0.6703 | 0.6692 | 0.6686 | 0.6682 |

2. Numerical experiment 2

We perform the second set of numerical experiment by adjusting boundary conditions of wave propagation in a rectangular beam from four free edges (FFFF) to two opposite free edges and two opposite clamped edges (FCFC), and three clamped edges and one free edge (CCFC) successively. The numerical results are presented in Tables 7 and 8. A brief discussion is as follows:

a. The dispersion curves of wave propagation are uniformly convergent. When the number of truncated series terms equals to 4, the computational errors are usually between 1.0E-3 and 4.0E-2, and show a trend of sustained decrease. When the number of truncated terms equals to 16, the computational errors further decrease by 1-2 orders of magnitude, and are usually between 1.0E-4 and 5.0E-3.

b. The change of boundary conditions has little effect on the uniform convergence of the dispersion curves of wave propagation.

Table 7: Computational error $e(\Omega^l_{M,N})$ of dispersion curve corresponding to FCFC boundaries.

| Wave | Branch Order $l$ | $M=N=4$ | $M=N=8$ | $M=N=12$ | $M=N=16$ |
|---|---|---|---|---|---|
| $L$ | 1 | 0.0062 | 0.0049 | 0.0045 | 0.0044 |
|  | 2 | 0.0028 | 0.0008 | 0.0003 | 0.0001 |
|  | 3 | 0.0036 | 0.0006 | 0.0002 | 0.0001 |
|  | 4 | 0.0109 | 0.0025 | 0.0010 | 0.0004 |
|  | 5 | 0.0088 | 0.0021 | 0.0008 | 0.0003 |
|  | 6 | 0.0110 | 0.0034 | 0.0014 | 0.0005 |
|  | 7 | 0.0074 | 0.0023 | 0.0009 | 0.0004 |
|  | 8 | 0.0017 | 0.0003 | 0.0001 | 0.0001 |
|  | 9 | 0.0063 | 0.0022 | 0.0009 | 0.0003 |
|  | 10 | 0.0161 | 0.0020 | 0.0006 | 0.0002 |
|  | 11 | 0.0275 | 0.0060 | 0.0023 | 0.0008 |
|  | 12 | 0.0323 | 0.0024 | 0.0009 | 0.0003 |
|  | 13 | 0.0323 | 0.0006 | 0.0003 | 0.0001 |
|  | 14 | 0.0017 | 0.0002 | 0.0001 | 0.0001 |
| $T$ | 1 | 0.0011 | 0.0002 | 0.0001 | 0.0001 |
|  | 2 | 0.0384 | 0.0126 | 0.0051 | 0.0018 |
|  | 3 | 0.0167 | 0.0057 | 0.0024 | 0.0009 |
|  | 4 | 0.0162 | 0.0037 | 0.0013 | 0.0005 |
|  | 5 | 0.0222 | 0.0055 | 0.0022 | 0.0008 |
|  | 6 | 0.0073 | 0.0025 | 0.0017 | 0.0008 |
|  | 7 | 0.0097 | 0.0030 | 0.0012 | 0.0004 |
|  | 8 | 0.0014 | 0.0012 | 0.0007 | 0.0003 |
|  | 9 | 0.0075 | 0.0051 | 0.0024 | 0.0009 |
|  | 10 | 0.0124 | 0.0035 | 0.0014 | 0.0005 |
|  | 11 | 0.0391 | 0.0157 | 0.0060 | 0.0021 |
|  | 12 | 0.0282 | 0.0091 | 0.0036 | 0.0013 |
|  | 13 | 0.0116 | 0.0003 | 0.0012 | 0.0006 |
| $B_{x1}$ | 1 | 0.0193 | 0.0071 | 0.0030 | 0.0011 |
|  | 2 | 0.0090 | 0.0038 | 0.0016 | 0.0007 |
|  | 3 | 0.0204 | 0.0047 | 0.0018 | 0.0007 |
|  | 4 | 0.0234 | 0.0061 | 0.0025 | 0.0009 |
|  | 5 | 0.0190 | 0.0060 | 0.0026 | 0.0009 |
|  | 6 | 0.0203 | 0.0105 | 0.0049 | 0.0018 |
|  | 7 | 0.0366 | 0.0047 | 0.0028 | 0.0027 |
|  | 8 | 0.0339 | 0.0058 | 0.0021 | 0.0008 |
|  | 9 | 0.0018 | 0.0002 | 0.0002 | 0.0001 |
|  | 10 | 0.0332 | 0.0067 | 0.0026 | 0.0009 |
|  | 11 | 0.0135 | 0.0025 | 0.0007 | 0.0003 |
|  | 12 | 0.0397 | 0.0214 | 0.0079 | 0.0028 |
|  | 13 | 0.0113 | 0.0127 | 0.0073 | 0.0030 |
| $B_{x2}$ | 1 | 0.0112 | 0.0028 | 0.0010 | 0.0004 |
|  | 2 | 0.0046 | 0.0013 | 0.0006 | 0.0002 |
|  | 3 | 0.0394 | 0.0062 | 0.0019 | 0.0006 |
|  | 4 | 0.0285 | 0.0079 | 0.0022 | 0.0007 |
|  | 5 | 0.0234 | 0.0037 | 0.0012 | 0.0004 |
|  | 6 | 0.0160 | 0.0027 | 0.0009 | 0.0003 |
|  | 7 | 0.0179 | 0.0055 | 0.0023 | 0.0009 |
| | 8 | 0.0398 | 0.0064 | 0.0022 | 0.0007 |
| | 9 | 0.0394 | 0.0100 | 0.0022 | 0.0006 |
| | 10 | 0.0396 | 0.0254 | 0.0074 | 0.0023 |
| | 11 | 0.0077 | 0.0018 | 0.0006 | 0.0002 |
| | 12 | 0.0175 | 0.0024 | 0.0012 | 0.0004 |
| | 13 | 0.0080 | 0.0029 | 0.0012 | 0.0005 |
| | 14 | 0.0204 | 0.0016 | 0.0004 | 0.0001 |
| | 15 | 0.0398 | 0.0050 | 0.0020 | 0.0007 |

Table 8: Computational error $e(\Omega_{M,N}^l)$ of dispersion curve corresponding to CCFC boundaries.

| Wave | Branch Order $l$ | $M=N=4$ | $M=N=8$ | $M=N=12$ | $M=N=16$ |
|---|---|---|---|---|---|
| L | 1 | 0.0040 | 0.0010 | 0.0020 | 0.0017 |
| | 2 | 0.0038 | 0.0011 | 0.0006 | 0.0005 |
| | 3 | 0.0060 | 0.0036 | 0.0031 | 0.0027 |
| | 4 | 0.0019 | 0.0006 | 0.0003 | 0.0002 |
| | 5 | 0.0161 | 0.0050 | 0.0032 | 0.0030 |
| | 6 | 0.0062 | 0.0015 | 0.0005 | 0.0005 |
| | 7 | 0.0061 | 0.0014 | 0.0005 | 0.0004 |
| | 8 | 0.0118 | 0.0040 | 0.0015 | 0.0008 |
| | 9 | 0.0076 | 0.0022 | 0.0011 | 0.0005 |
| | 10 | 0.0043 | 0.0012 | 0.0006 | 0.0003 |
| | 11 | 0.0044 | 0.0014 | 0.0006 | 0.0004 |
| | 12 | 0.0047 | 0.0020 | 0.0013 | 0.0009 |
| | 13 | 0.0075 | 0.0028 | 0.0013 | 0.0005 |
| | 14 | 0.0032 | 0.0017 | 0.0006 | 0.0003 |
| | 15 | 0.0193 | 0.0045 | 0.0015 | 0.0007 |
| | 16 | 0.0043 | 0.0015 | 0.0004 | 0.0003 |
| | 17 | 0.0277 | 0.0069 | 0.0025 | 0.0009 |
| | 18 | 0.0326 | 0.0075 | 0.0027 | 0.0009 |
| | 19 | 0.0379 | 0.0062 | 0.0021 | 0.0008 |
| | 20 | 0.0329 | 0.0024 | 0.0008 | 0.0003 |
| | 21 | 0.0067 | 0.0015 | 0.0006 | 0.0003 |
| | 22 | 0.0109 | 0.0022 | 0.0008 | 0.0004 |
| | 23 | 0.0086 | 0.0020 | 0.0009 | 0.0006 |
| | 24 | 0.0085 | 0.0020 | 0.0008 | 0.0005 |
| | 25 | 0.0104 | 0.0028 | 0.0024 | 0.0003 |
| $B_{x1}$ | 1 | 0.0134 | 0.0052 | 0.0022 | 0.0008 |
| | 2 | 0.0019 | 0.0010 | 0.0006 | 0.0004 |
| | 3 | 0.0135 | 0.0046 | 0.0019 | 0.0007 |
| | 4 | 0.0134 | 0.0040 | 0.0017 | 0.0006 |
| | 5 | 0.0100 | 0.0031 | 0.0012 | 0.0006 |
| | 6 | 0.0107 | 0.0067 | 0.0040 | 0.0031 |
| | 7 | 0.0050 | 0.0014 | 0.0007 | 0.0005 |
| | 8 | 0.0060 | 0.0014 | 0.0004 | 0.0003 |
| | 9 | 0.0040 | 0.0026 | 0.0014 | 0.0006 |
| | 10 | 0.0036 | 0.0036 | 0.0022 | 0.0010 |

| | | | | |
|---|---|---|---|---|
| 11 | 0.0027 | 0.0159 | 0.0158 | 0.0006 |
| 12 | 0.0014 | 0.0003 | 0.0001 | 0.0002 |
| 13 | 0.0031 | 0.0007 | 0.0007 | 0.0004 |
| 14 | 0.0046 | 0.0010 | 0.0005 | 0.0003 |
| 15 | 0.0106 | 0.0030 | 0.0013 | 0.0005 |
| 16 | 0.0051 | 0.0018 | 0.0009 | 0.0004 |
| 17 | 0.0009 | 0.0002 | 0.0002 | 0.0001 |
| 18 | 0.0065 | 0.0032 | 0.0014 | 0.0007 |
| 19 | 0.0119 | 0.0038 | 0.0016 | 0.0007 |
| 20 | 0.0048 | 0.0004 | 0.0002 | 0.0002 |
| 21 | 0.0030 | 0.0017 | 0.0007 | 0.0003 |
| 22 | 0.0336 | 0.0064 | 0.0019 | 0.0010 |
| 23 | 0.0140 | 0.0087 | 0.0036 | 0.0017 |

*5.3. Brief conclusions*

For convergence characteristics of the Fourier series multiscale solution of wave propagation in a beam of rectangular cross section, we arrive at concluding remarks as follows:

a. The Fourier series multiscale solution has rapid convergence speed, gives accurate results not only for lower order branches of dispersion curves, but also for higher order branches of dispersion curves.

b. The changes of length width ratio and boundary conditions have little effect on the convergence of the Fourier series multiscale solution.

## 6. Propagation characteristics of waves in a square beam

In this section, the inherent characteristics of wave propagation, such as frequency spectrum and wave modes, in the beam of square cross section, are investigated by the Fourier series multiscale solution with the number of truncated terms $M = N = 20$.

*6.1. Frequency spectrum*

As to the three types of boundary conditions, i.e., four free edges (FFFF), two opposite free edges and two opposite clamped edges (FCFC), and three clamped edges and one free edge (CCFC), wave propagation in the beam of square cross section is taken into consideration, and the frequency spectrums are presented in Figures 3-5. It is a pity that, the Fourier series multiscale method cannot accurately recognize the third order dispersion curve ($L_a^3$) and the fourth order dispersion curve ($L_a^4$) of $L_a$ waves in Figure 3. Accordingly, blind area appears near the position of intersection, $K \in [0.43, 0.60]$, of two dispersion curves.

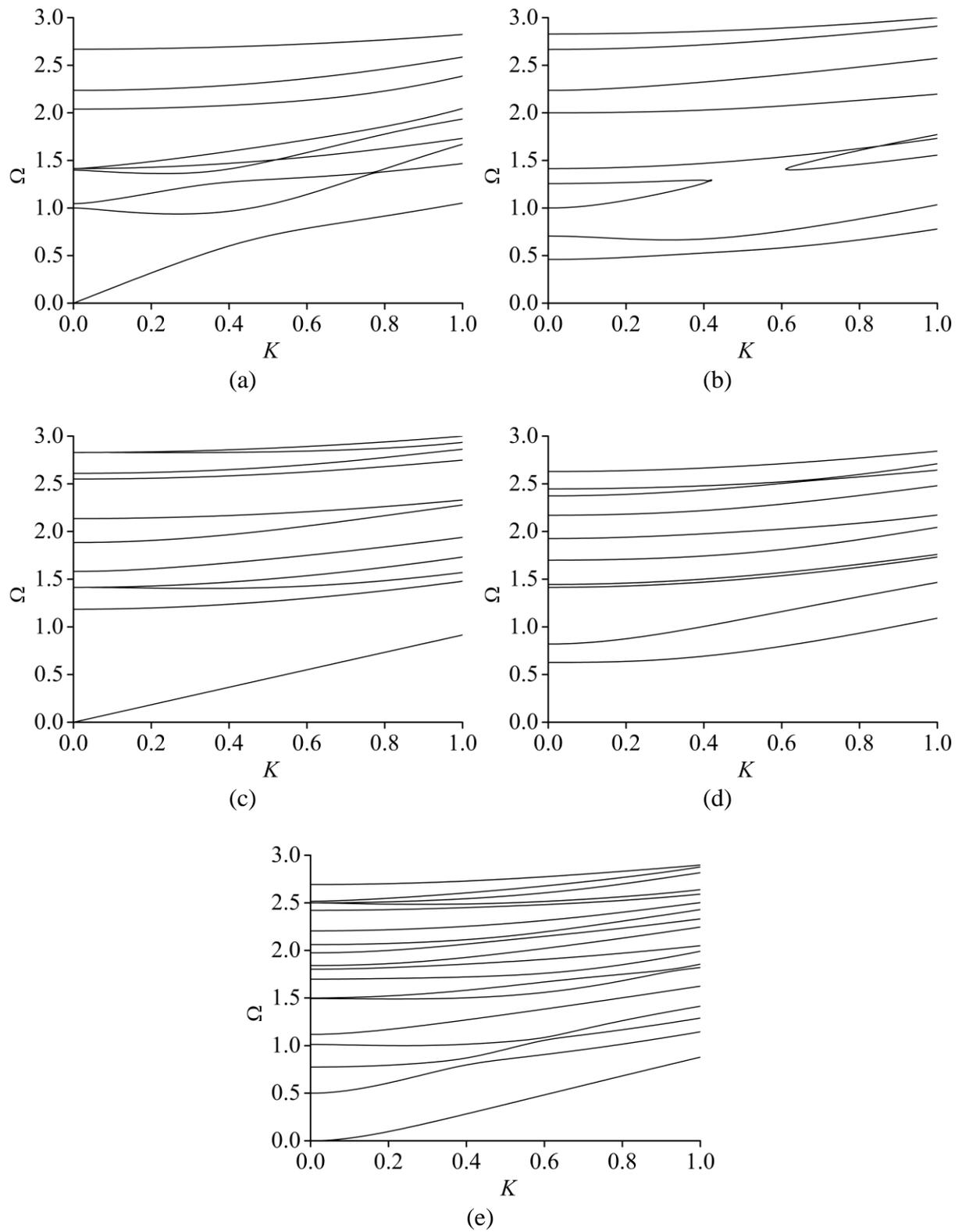

Figure 3: Frequency spectra corresponding to FFFF boundaries:
(a) $L_s$ waves, (b) $L_a$ waves, (c) $T_a$ waves, (d) $T_s$ waves, (e) $B_{x1}$ waves.

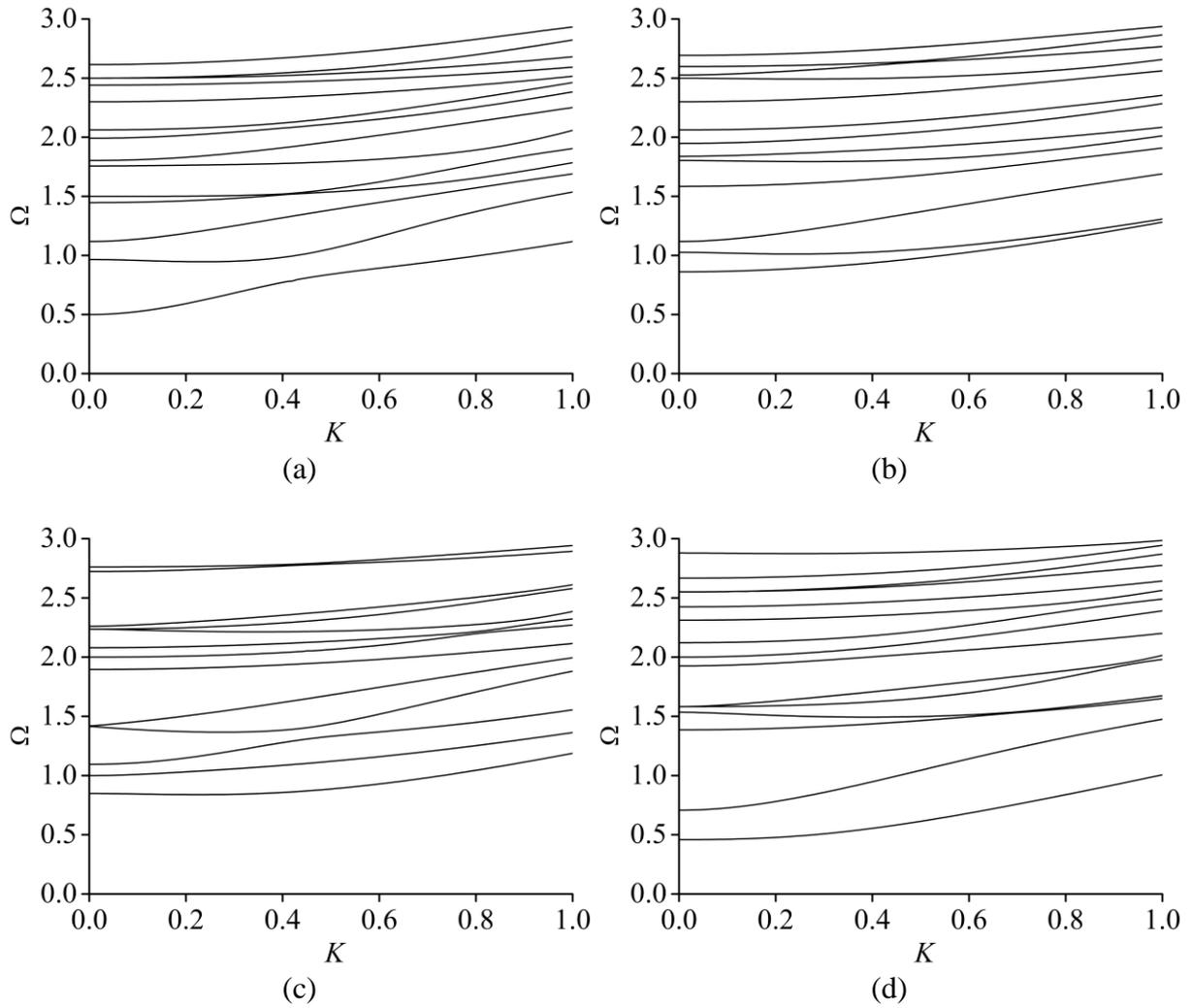

Figure 4: Frequency spectra corresponding to FCFC boundaries:
(a) $L$ waves, (b) $T$ waves, (c) $B_{x1}$ waves, (d) $B_{x2}$ waves.

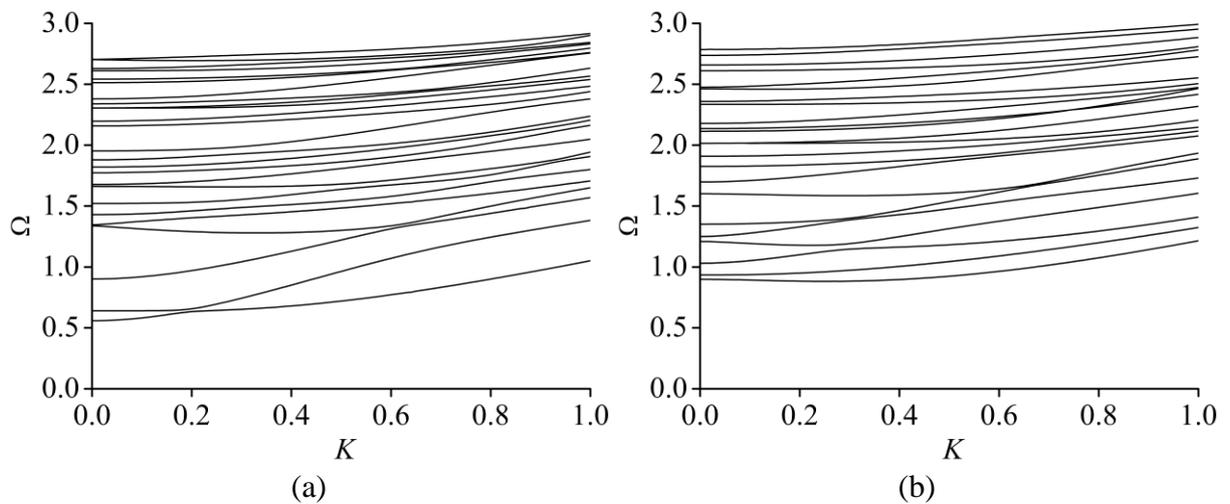

Figure 5: Frequency spectra corresponding to CCFC boundaries:
(a) $L$ waves, (b) $B_{x1}$ waves.

## 6.2. Wave modes

As to the dispersion curves of different types of traveling waves in the frequency spectrums, wave modes, namely deformation patterns of the cross section, corresponding to some characteristic points on the dispersion curves are given in Tables 9-17.

Table 9: Distortion and out-of-plane displacement corresponding to characteristic points on dispersion curves of $L_s$ and $L_a$ waves in FFFF boundaries.

| Points | Distortion | out-of-plane displacement | Points | Distortion | out-of-plane displacement |
|---|---|---|---|---|---|
| $L_s^1$, $K=0.1$ | | | $L_s^1$, $K=0.9$ | | |
| $L_s^2$, $K=0.1$ | | | $L_s^2$, $K=0.9$ | | |
| $L_s^3$, $K=0.1$ | | | $L_s^3$, $K=0.9$ | | |
| $L_s^9$, $K=0.1$ | | | $L_s^9$, $K=0.9$ | | |
| $L_a^1$, $K=0.1$ | | | $L_a^1$, $K=0.9$ | | |
| $L_a^2$, $K=0.1$ | | | $L_a^2$, $K=0.9$ | | |
| $L_a^3$, $K=0.1$ | | | $L_a^3$, $K=0.9$ | | |

Table 10: Distortion and out-of-plane displacement corresponding to characteristic points on dispersion curves of $T_a$ and $T_s$ waves in FFFF boundaries.

| Points | Distortion | out-of-plane displacement | Points | Distortion | out-of-plane displacement |
|---|---|---|---|---|---|
| $T_a^1$, $K = 0.1$ | | | $T_a^1$, $K = 0.9$ | | |
| $T_a^2$, $K = 0.1$ | | | $T_a^2$, $K = 0.9$ | | |
| $T_a^3$, $K = 0.1$ | | | $T_a^3$, $K = 0.9$ | | |
| $T_a^5$, $K = 0.1$ | | | $T_a^5$, $K = 0.9$ | | |
| $T_s^1$, $K = 0.1$ | | | $T_s^1$, $K = 0.9$ | | |
| $T_s^2$, $K = 0.1$ | | | $T_s^2$, $K = 0.9$ | | |
| $T_s^5$, $K = 0.1$ | | | $T_s^5$, $K = 0.9$ | | |
| $T_s^{10}$, $K = 0.1$ | | | $T_s^{10}$, $K = 0.9$ | | |

Table 11: Distortion and out-of-plane displacement corresponding to characteristic points on dispersion curves of $B_{x1}$ waves in FFFF boundaries.

| Points | Distortion | out-of-plane displacement | Points | Distortion | out-of-plane displacement |
|---|---|---|---|---|---|
| $B_{x_1}^1$ $K=0.1$ | | | $B_{x_1}^1$ $K=0.9$ | | |
| $B_{x_1}^2$ $K=0.1$ | | | $B_{x_1}^2$ $K=0.9$ | | |
| $B_{x_1}^3$ $K=0.1$ | | | $B_{x_1}^3$ $K=0.9$ | | |
| $B_{x_1}^4$ $K=0.1$ | | | $B_{x_1}^4$ $K=0.9$ | | |
| $B_{x_1}^5$ $K=0.1$ | | | $B_{x_1}^5$ $K=0.9$ | | |
| $B_{x_1}^{18}$ $K=0.1$ | | | $B_{x_1}^{18}$ $K=0.9$ | | |

Table 12: Distortion and out-of-plane displacement corresponding to characteristic points on dispersion curves of *L* waves in FCFC boundaries.

| Points | Distortion | out-of-plane displacement | Points | Distortion | out-of-plane displacement |
|---|---|---|---|---|---|
| $L^1$ $K=0.1$ | | | $L^1$ $K=0.9$ | | |
| $L^2$ $K=0.1$ | | | $L^2$ $K=0.9$ | | |
| $L^3$ $K=0.1$ | | | $L^3$ $K=0.9$ | | |
| $L^4$ $K=0.1$ | | | $L^4$ $K=0.9$ | | |
| $L^5$ $K=0.1$ | | | $L^5$ $K=0.9$ | | |
| $L^{14}$ $K=0.1$ | | | $L^{14}$ $K=0.9$ | | |

Table 13: Distortion and out-of-plane displacement corresponding to characteristic points on dispersion curves of *T* waves in FCFC boundaries.

| Points | Distortion | out-of-plane displacement | Points | Distortion | out-of-plane displacement |
|---|---|---|---|---|---|
| $T^1$ $K=0.1$ | 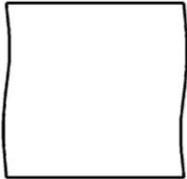 | 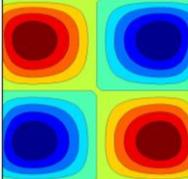 | $T^1$ $K=0.9$ | 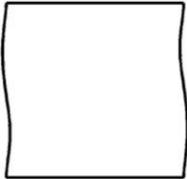 | 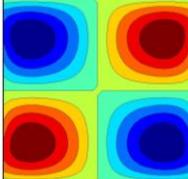 |
| $T^2$ $K=0.1$ | 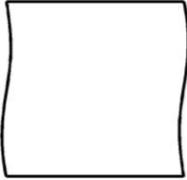 | 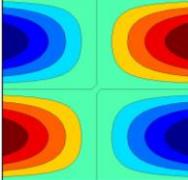 | $T^2$ $K=0.9$ | 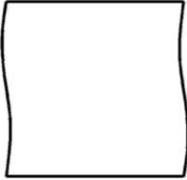 | 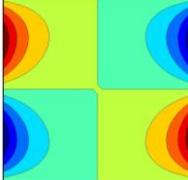 |
| $T^3$ $K=0.1$ | 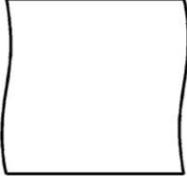 | 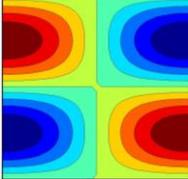 | $T^3$ $K=0.9$ | 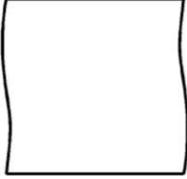 | 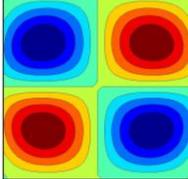 |
| $T^4$ $K=0.1$ | 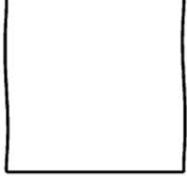 | 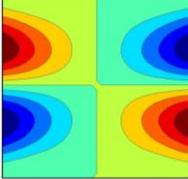 | $T^4$ $K=0.9$ | 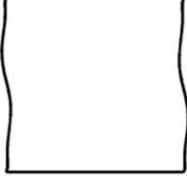 | 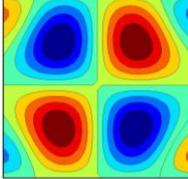 |
| $T^5$ $K=0.1$ | 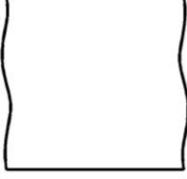 | 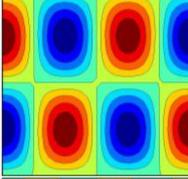 | $T^5$ $K=0.9$ | 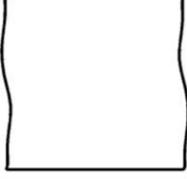 | 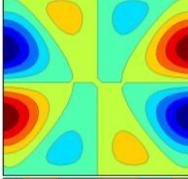 |
| $T^{13}$ $K=0.1$ | 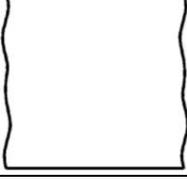 | 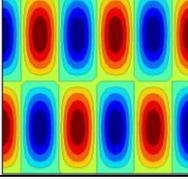 | $T^{13}$ $K=0.9$ | 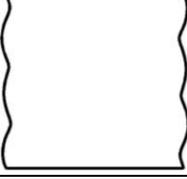 | 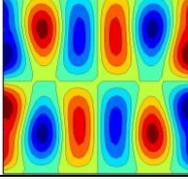 |

Table 14: Distortion and out-of-plane displacement corresponding to characteristic points on dispersion curves of $B_{x1}$ waves in FCFC boundaries.

| Points | Distortion | out-of-plane displacement | Points | Distortion | out-of-plane displacement |
|---|---|---|---|---|---|
| $B_{x_1}^{1}$ $K=0.1$ | | | $B_{x_1}^{1}$ $K=0.9$ | | |
| $B_{x_1}^{2}$ $K=0.1$ | | | $B_{x_1}^{2}$ $K=0.9$ | | |
| $B_{x_1}^{3}$ $K=0.1$ | | | $B_{x_1}^{3}$ $K=0.9$ | | |
| $B_{x_1}^{4}$ $K=0.1$ | | | $B_{x_1}^{4}$ $K=0.9$ | | |
| $B_{x_1}^{5}$ $K=0.1$ | | | $B_{x_1}^{5}$ $K=0.9$ | | |
| $B_{x_1}^{13}$ $K=0.1$ | | | $B_{x_1}^{13}$ $K=0.9$ | | |

Table 15: Distortion and out-of-plane displacement corresponding to characteristic points on dispersion curves of $B_{x2}$ waves in FCFC boundaries.

| Points | Distortion | out-of-plane displacement | Points | Distortion | out-of-plane displacement |
|---|---|---|---|---|---|
| $B_{x_2}^1$ $K=0.1$ | | | $B_{x_2}^1$ $K=0.9$ | | |
| $B_{x_2}^2$ $K=0.1$ | | | $B_{x_2}^2$ $K=0.9$ | | |
| $B_{x_2}^3$ $K=0.1$ | | | $B_{x_2}^3$ $K=0.1$ | | |
| $B_{x_2}^4$ $K=0.1$ | | | $B_{x_2}^4$ $K=0.9$ | | |
| $B_{x_2}^5$ $K=0.1$ | | | $B_{x_2}^5$ $K=0.9$ | | |
| $B_{x_2}^{13}$ $K=0.1$ | | | $B_{x_2}^{13}$ $K=0.9$ | | |

Table 16: Distortion and out-of-plane displacement corresponding to characteristic points on dispersion curves of *L* waves in CCFC boundaries.

| Points | Distortion | out-of-plane displacement | Points | Distortion | out-of-plane displacement |
|---|---|---|---|---|---|
| $L^1$ $K=0.1$ | | | $L^1$ $K=0.9$ | | |
| $L^2$ $K=0.1$ | | | $L^2$ $K=0.9$ | | |
| $L^3$ $K=0.1$ | | | $L^3$ $K=0.9$ | | |
| $L^4$ $K=0.1$ | | | $L^4$ $K=0.9$ | | |
| $L^5$ $K=0.1$ | | | $L^5$ $K=0.9$ | | |
| $L^{25}$ $K=0.1$ | | | $L^{25}$ $K=0.9$ | | |

Table 17: Distortion and out-of-plane displacement corresponding to characteristic points on dispersion curves of $B_{x1}$ waves in CCFC boundaries.

| Points | Distortion | out-of-plane displacement | Points | Distortion | out-of-plane displacement |
|---|---|---|---|---|---|
| $B_{x_1}^1$ $K=0.1$ | | | $B_{x_1}^1$ $K=0.9$ | | |
| $B_{x_1}^2$ $K=0.1$ | | | $B_{x_1}^2$ $K=0.9$ | | |
| $B_{x_1}^3$ $K=0.1$ | | | $B_{x_1}^3$ $K=0.9$ | | |
| $B_{x_1}^4$ $K=0.1$ | | | $B_{x_1}^4$ $K=0.9$ | | |
| $B_{x_1}^5$ $K=0.1$ | | | $B_{x_1}^5$ $K=0.9$ | | |
| $B_{x_1}^{23}$ $K=0.1$ | | | $B_{x_1}^{23}$ $K=0.9$ | | |

Now we give a brief analysis of the evolution of wave modes:

a. For the small values of wave number, the wave modes corresponding to the first branches of the dispersion curves, namely $L_s^1$, $T_a^1$ and $B_{x_1}^1$, of the three kinds of waves in the square beam with FFFF boundary conditions satisfy the hypothesis of plane sections well.

b. With the values of wave number increasing, the wave modes corresponding to the lower order branches of dispersion curves tend to be complicated. Meanwhile, for small or large values of wave number, the wave modes corresponding to the higher order branches of dispersion curves are complicated. And in some specific cases, sharp gradients or singularity appears as the typical multiscale phenomenon.

The usually used one-dimensional simplified models, such as the longitudinal vibration model of rods, the torsional vibration model of rods, and the transverse vibration model of straight beams, are all based on the hypothesis of plane sections. Accordingly, the foregoing analysis of wave modes verifies that these one-dimensional simplified models are only suitable for use in the approximate analysis of the first order traveling waves in square beam subjected to FFFF boundary conditions when the wave number takes small values.

# 7. Conclusions

The problem of elastic wave propagation constitutes an extremely important area of study in geophysics, crystal physics, electronics, vibration control, and nondestructive testing. In this paper, the exact analysis of wave propagation in a rectangular beam is extended to a thorough multiscale analysis for a system of completely coupled second order linear differential equations for modal functions, where general boundary conditions are prescribed. It is concluded that:

1. We derive the Fourier series multiscale solution for wave propagation in the rectangular beam, by means of which an accurate wave propagation model is established.

2. We discover the inherent characteristics, such as frequency spectrum, cut-off frequencies and wave modes, of wave propagation in the square beam, which enrich the research of the classical elastodynamics.

3. We analyze the multiscale characteristics of wave modes in the square beam and redefine the scope of application of the usually simplified one-dimensional models.

The preliminary study on application to exact analysis of wave propagation in a beam with rectangular cross section verifies the effectiveness of the present Fourier series multiscale method, and the accurate wave propagation model established potentially offers a powerful means for the simultaneous control of coupled waves in the beam, as well as the development of guided wave NDE techniques.